%% file: draft.tex
\documentclass[10pt]{article}
\usepackage[letterpaper,twoside,outer=1.3in,vmargin=1.3in,]{geometry}

\usepackage{amsmath,amsfonts,amsthm,url,graphicx,subcaption,float,booktabs,bbm}
\usepackage[ruled,vlined,linesnumbered]{algorithm2e}

\def \bR {\mathbb{R}}
\def \bZ {\mathbb{Z}}
\def \bL {\mathcal{L}}
\def \conv {\text{conv}}
\def \proj {\text{proj}}
\def \bQs {\overline{Q}_s^*}
\def \hQs {Q_{s}^*}
\def \THE {\text{Restricted Separation of Lagrangian Cuts}}
\def \strben {\emph{StrBen}}
\def \exact {\emph{Exact}}
\def \rstrone {\emph{Rstr1}}
\def \rstrtwo {\emph{Rstr2}}
\def \rstrmip {\emph{RstrMIP}}
\def \BDD {\emph{BDD}}
\def \xslambda {\bar{x}^s_\lambda}
\def \yslambda {\bar{y}^s_\lambda}
\def \Ps#1{P_s(#1)}
\def \PsPis {P_s(\Pi_s)}
\def \PsPistar {P_s(\Pi_s^*)}

\newcommand{\gv}[1]{#1}
\newcommand{\pv}[1]{}

\newtheorem{thm}{Theorem}

\newtheorem{prop}[thm]{Proposition}

\theoremstyle{definition}

\author{Rui Chen, James Luedtke}
\title{On Generating Lagrangian Cuts for Two-Stage Stochastic Integer Programs}
\date{\small{Department of Industrial and Systems Engineering, University of Wisconsin-Madison}\\
}

\begin{document}
	\maketitle
	\begin{abstract}
		We investigate new methods for generating Lagrangian cuts to solve two-stage stochastic integer programs. Lagrangian
		cuts can be added to a Benders reformulation, and are derived from solving single scenario integer
		programming subproblems identical to those used in the nonanticipative Lagrangian dual of a stochastic integer program. 
		While Lagrangian cuts have the potential to significantly strengthen the Benders relaxation,
		generating Lagrangian cuts can be computationally demanding. We investigate new techniques for generating Lagrangian
		cuts with the goal of obtaining methods that provide significant improvements to the Benders relaxation quickly.
		Computational results demonstrate that our proposed method improves the Benders relaxation  significantly faster than
		previous methods for generating Lagrangian cuts and, when used within a branch-and-cut algorithm,  significantly reduces
		the size of the search tree for three classes of test problems.
	\end{abstract}
	
	\input{joc/main}
	\appendix
	\section*{Appendix}
	\subsection*{A Branch-and-Cut Algorithm}\label{BnC}
	A standard branch-and-cut algorithm for solving two-stage SIPs with either continuous recourse or pure binary
	first-stage variables is described in Algorithm \ref{alg:BnC}.
	The algorithm starts with the approximations $\{\hat{Q}_s\}_{s\in S}$ of $\{Q_s\}_{s\in S}$ uses it to create an LP
	initial LP relaxation of problem \eqref{Benders} represented by the root node $0$ which is added to the list of
	candidate nodes $\mathcal{L}$. In our case, the approximations
	$\{\hat{Q}_s\}_{s\in S}$ are constructed by Benders cuts and Lagrangian cuts. Within each iteration, the algorithm picks
	a node from  $\mathcal{L}$ and solves the corresponding LP which is a relaxation of problem \eqref{Benders} with the addition of the
	constraints adding from branching,  and generates a candidate solution $(\hat{x},\{\hat{\theta}_s\}_{s\in S})$. If
	$\hat{x}$ does not satisfy the integrality constraints ($\hat{x} \notin X$), the algorithm creates two new subproblems (nodes)
	by subdividing the feasible solutions to that node. The subproblems are created by choosing an
	integer decision variable $x_j$ having $\hat{x}_j \notin \mathbb{Z}$, and adding the constraint $x_j \leq \lfloor
	\hat{x}_j \rfloor$ in one of the subproblems and adding the constraint $x_j \geq \lceil \hat{x}_j \rceil$ in the other
	subproblem.  If $\hat{x}$ satisfies the integrality constraints ($\hat{x} \in X$), the algorithm evaluates $Q_s(\hat{x})$ for all $s\in S$.
	If any violated Benders cuts and integer L-shaped cuts are identified, these are used to update $\{\hat{Q}_s\}_{s\in
		S}$, and the relaxation
	at that node is then solved again. This process of searching for and adding cuts if violated when $\hat{x}
	\in X$ can be implemented using the ``Lazy constraint callback'' in commercial solvers.  If no violated Benders or integer L-shaped cuts are identified, then the current
	solution is feasible and thus the upper bound $\bar{z}$ and best known solution (the incumbent $x^*$) is updated. The algorithm
	terminates when there are no more nodes to explore in the candidate list $\mathcal{L}$.
	\begin{algorithm}[H]\label{alg:BnC}
		\SetAlgoLined
		\KwIn{$\{\hat{Q}_s\}_{s\in S}$}
		\KwOut{$(x^*,\{\theta_s^* \}_{s\in S})$}
		Initialize $\bL\leftarrow\{0\}$,  LP$_0\leftarrow\min_{x,\theta_s}\{c^Tx+\sum_{s\in S}p_s\theta_s:\theta_s\geq\hat{Q}_s(x),s\in S,Ax\geq b\}$, $\bar{z}\leftarrow+\infty$, $(x^*,\{\theta_s^* \}_{s\in S})\leftarrow\emptyset$\\
		\While{$\bL\neq\emptyset$}{
			Choose a node $i\in\bL$\\
			Solve LP$_i$ and, if feasible, obtain optimal solution $(\hat{x},\{\hat{\theta}_s\}_{s\in S})$ and optimal value $\hat{z}$\\
			\eIf{LP$_i$ is feasible and $\hat{z}<\bar{z}$}{
				\eIf{$\hat{x}\in X$}{
					\For{$s\in S$}{
						Evaluate $Q_s(\hat{x})$\\
						Add the Benders cut \eqref{BenCut} to update $\hat{Q}_s$ if violated\\
						\If{$X=\{0,1\}^n$}{Add the integer L-shaped cut \eqref{IntLShaped} to update $\hat{Q}_s$ if violated}
					}
					\If{$\hat{\theta}_s=Q_s(\hat{x})$ for all $s\in S$}{
						$\bL\leftarrow\bL\setminus\{i\}$\\
						$(x^*,\{\theta_s^* \}_{s\in S})\leftarrow(\hat{x},\{\hat{\theta}_s\}_{s\in S})$\\
						$\bar{z}\leftarrow\hat{z} $}
				}
				{Choose an integer variable $x_j$ with $\hat{x}_j\neq\bZ$, then branch on LP$_i$ to construct LP$_{i_1}$ and LP$_{i_2}$ such that LP$_{i_1}$ is LP$_i$ with an additional constraint $x_j\leq\lfloor \hat{x}_j\rfloor$ and LP$_{i_2}$ is LP$_i$ with an additional constraint $x_j\geq\lceil \hat{x}_j\rceil$\\
					$\bL\leftarrow (\bL\setminus\{i\})\cup\{i_1,i_2 \}$}
			}
			{$\bL\leftarrow\bL\setminus\{i\}$}
		}
		\caption{Branch-and-cut for SIPs.}
	\end{algorithm}

	\subsection*{Test Problems}\label{test_prob}
	\subsubsection*{The SSLP problem}
	The SSLP problem \cite{ntaimo2005million} is a two-stage SIP with pure binary first-stage and mixed-binary second-stage variables.  
	In this problem, the decision maker has to choose from $m$ sites to allocate servers with some cost in the first stage. Then in the second stage, the availability of each client $i$ would be observed and every available client must be served at some site. The first-stage variables are denoted by $x$ with $x_j=1$ if and only if a server is located at site $j$. The second-stage variables are denoted by $y$ and $y_0$ where $y_{ij}^s=1$ if and only if client $i$ is served at site $j$ in scenario $s$ and $y_{0j}^s$ is the amount of resource shortage at at site $j$ in scenario $s$. Allocating a server costs $c_j$ at site $j$. Each allocated server can provide up to $u$ units of resource. Client $i$ uses $d_{ij}$ units of resource and generates revenue $q_{ij}$ if served at site $j$. Each unit of resource shortage at site $j$ incurs a penalty $q_{0j}$. The client availability in scenario $s$ is represented by $h^s$ with $h^s_i=1$ if and only if client $i$ is present in scenario $s$. The extensive formulation of the problem is as follows:\begin{displaymath}
		\begin{aligned}
			\min_{x,y^s,y_0^s}\ & \sum_{j=1}^mc_jx_j+\sum_{s\in S}p_s\Big(\sum_{j=1}^m q_{0j}y_{0j}^s-\sum_{i=1}^n\sum_{j=1}^mq_{ij}y_{ij}^s\Big)\\
			\text{s.t. }\ &\sum_{i=1}^nd_{ij}y_{ij}^s-y_{0j}^s\leq ux_j, &&j=1,\ldots,m,\ s\in S,\\
			&\sum_{j=1}^my_{ij}^s=h_i^s, &&i=1,\ldots,n,\ s\in S,\\
			&x_{j}\in\{0,1\}, &&j=1,\ldots,m,\\
			&y_{ij}^s\in\{0,1\}, &&i=1,\ldots,n,\ j=1,\ldots,m,\ s\in S,\\
			&y_{0j}^s\geq 0, &&j=1,\ldots,m,\ s\in S.
		\end{aligned}
	\end{displaymath}
	All SSLP instances are generated with $(m,n)\in \{(20,100),(30,70),(40,50),(50,40)\}$ and $|S|\in \{50,200\}$. For each size $(m,n,|S|)$, three instances are generated with $k\in\{1,2,3\}$. For each SSLP instance, the data are generated as:\begin{itemize}
		\item $p_s=1/|S|$, $s\in S$.
		\item $c_j\sim$Uniform($\{40,41,\ldots,80\}$), $j=1,\ldots,m$.
		\item $d_{ij}=q_{ij}\sim$Uniform($\{0,1,\ldots,25\}$), $i=1,\ldots,n$, $j=1,\ldots,m$.
		\item $q_{0j}=1000$, $j=1,\ldots,m$.
		\item $u=\sum_{i=1}^n\sum_{j=1}^m d_{ij}/m$.
		\item $h_i^s\sim$Bernoulli($1/2$), $i=1,\ldots,n$, $s\in S$.
	\end{itemize}
	All SSLP instances are named of the form sslp$k$\_$m$\_$n$\_$|S|$ where $m$ denotes the number of sites, $n$ denotes number of clients, $|S|$ is the number of scenarios and $k$ is the instance number.
	
	\subsubsection*{The SSLPV problem}
	The SSLPV problem is a variant of the SSLP problem that has mixed-binary first-stage and pure continuous second-stage
	variables. There are two main differences between SSLP and SSLPV. First, in the SSLPV problem, in addition to choosing
	which sites to ``open'' for servers, the decision maker also has an option to add more  servers than the ``base number''
	if the site is open. Second, in the second stage, we allow the clients' demands to be partly met by the servers, and
	only unmet shortage will be penalized. Thus, in this problem the recourse problem is a linear program. 
	We assume that if a site $j$ is selected (binary decision) it has base capacity of $u^{b}$ (independent of $j$), and installing this has cost
	$c^b_j$. In addition, it is possible to install any fraction of an additional $u^{n}$ units of capacity at site $j$, with total cost $c_j^n$ if the total new
	capacity is installed and the cost reduced proportionallly if a fraction is installed.
	The variable representing the amount of additional server capacity added at site $j$ is denoted by a continuous variable $x^+_j\in[0,1]$. The extensive formulation of the problem is as follows:\begin{displaymath}
		\begin{aligned}
			\min_{x,x^+,y^s,y_0^s}\ & \sum_{j=1}^m\Big(c_j^b x_j+ c_j^n x^+_j\Big)+\sum_{s\in S}p_s\Big(\sum_{j=1}^m q_{0j}y_{0j}^s-\sum_{i=1}^n\sum_{j=1}^mq_{ij}y_{ij}^s\Big)\\
			\text{s.t. }\ &\sum_{i=1}^nd_{ij}y_{ij}^s-y_{0j}^s\leq u^b x_j+ u^n x^+_j, &&j=1,\ldots,m,\ s\in S,\\
			&\sum_{j=1}^my_{ij}^s=h_i^s, &&i=1,\ldots,n,\ s\in S,\\
			&x_{j}\in\{0,1\},~x^+_j\in[0,1],~x^+_j\leq x_j, &&j=1,\ldots,m,\\
			&y_{ij}^s\in[0,1], &&i=1,\ldots,n,\ j=1,\ldots,m,\ s\in S,\\
			&y_{0j}^s\geq 0, &&j=1,\ldots,m,\ s\in S.
		\end{aligned}
	\end{displaymath}
	Each SSLPV instance corresponds to a SSLP instance.
	In all SSLPV instances, we set $u^b=(1-\rho)u$, $u^n = \rho u$, and for each site $j$, we set $c_j^b = (1-\rho)c_j$,
	$c_j^n = \rho c_j$, where $\rho=0.5$ and $u$ and $c_j$ are the data from the SSLP instance. The other data is the same as the data of the SSLP instances. All SSLP instances are named of the form sslp$k$\_$m$\_$n$\_$|S|$\_var where $m$ denotes the number of sites, $n$ denotes number of clients, $|S|$ is the number of scenarios and $k$ is the instance number.
	
	\subsubsection*{The SNIP problem}
	The SNIP problem \cite{pan2008minimizing} is a two-stage SIP with pure binary first-stage and continuous second-stage
	variables. In this problem the defender interdicts arcs on a directed network by installing sensors in the first stage
	to maximize the probability of finding an attacker. Then in the second stage, the origin and destination of
	the attacker  are observed and the attacker chooses the maximum reliability path from its origin to its destination. Let
	$N$ and $A$ denote the node set and the arc set of the network and let $D\subseteq A$ denote the set of interdictable
	arcs. The first-stage variables are denoted by $x$ with $x_a=1$ if and only if the defender installs a sensor on arc
	$a\in D$. The scenarios $s \in S$ define the possible origin/destination combinations of the attacker, with $u^s$
	representing the origin and $v^s$ representing the destination of the attacker for each scenario $s \in S$. The
	second-stage variables are denoted by $\pi$ where $\pi^s_i$ denotes the maximum probability of reaching destination
	$v^s$ undetected from node $i$ in scenario $s$. The budget for installing sensors is $b$, and the cost of installing a
	sensor on arc $a$ is $c_a$ for each arc $a\in D$. For each arc $a\in A$, the probability of traveling on arc $a$
	undetected is $r_{a}$ if the arc is not interdicted, or $q_a$ if the arc is interdicted. Finally, let $\bar{\pi}_j^s$ denote the maximum probability of reaching the destination undetected from node $j$ when no sensors are installed. The extensive formulation of the problem is as follows:\begin{displaymath}
		\begin{aligned}
			\min_{x,\pi^s}\ &\sum_{s\in S}p_s\pi_{u^s}^s\\
			\text{s.t. }\ &\sum_{a\in A}c_ax_a\leq b,\\
			&\pi_i^s-r_a\pi_j^s\geq 0, &&a=(i,j)\in A\setminus D,\ s\in S,\\
			&\pi_i^s-r_a\pi_j^s\geq -(r_a-q_a)\bar{\pi}_j^sx_a, &&a=(i,j)\in D,\ s\in S,\\
			&\pi_i^s-q_a\pi_j^s\geq 0, &&a=(i,j)\in D,\ s\in S,\\
			&\pi_{v^s}^s=1, &&s\in S,\\
			&x_a\in \{0,1\}, &&a\in D.
		\end{aligned}
	\end{displaymath}
	All of the SNIP instances we used in this paper are from \cite{pan2008minimizing}. We consider instances with snipno = 3, and budget $b\in\{30,50,70,90\}$. All instances have 320 first-stage binary variables, 2586 second-stage continuous variables per scenario and 456 scenarios.
	
	\subsection*{Some Convergence Profiles}\label{cvg_profiles}
	
\begin{figure}[H]
		\centering
		\begin{subfigure}[b]{0.38\linewidth}
			\includegraphics[width=\linewidth]{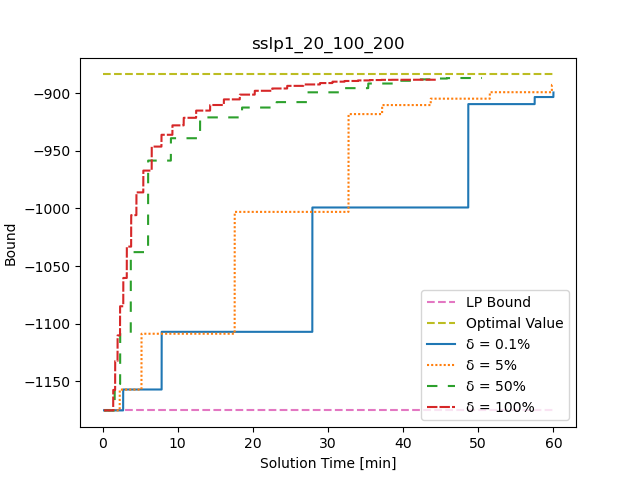}
		\end{subfigure}
		\begin{subfigure}[b]{0.38\linewidth}
			\includegraphics[width=\linewidth]{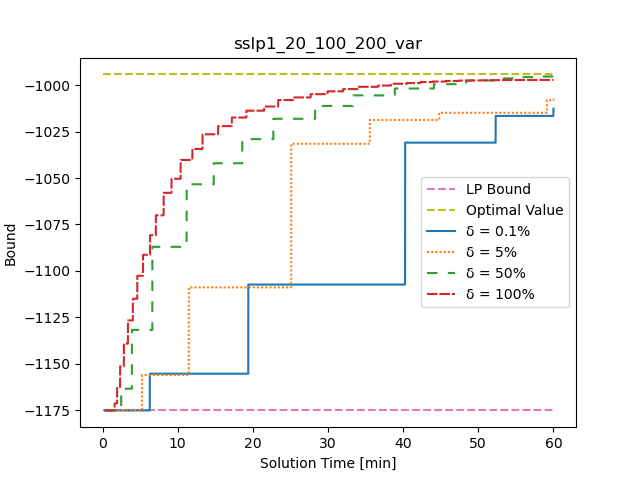}
		\end{subfigure}
		\begin{subfigure}[b]{0.38\linewidth}
			\includegraphics[width=\linewidth]{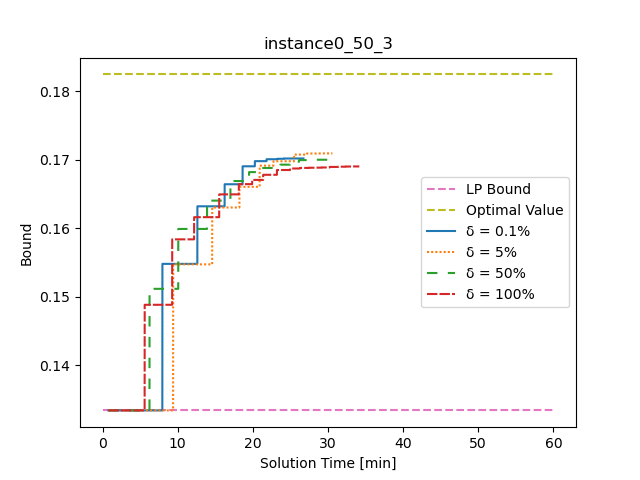}
		\end{subfigure}
		\caption{Convergence profile of {\rstrone} on an SSLP instance (top left), {an SSLPV instance (top right)}, and a
			SNIP instance (bottom) with $K=10$ and varying $\delta$ values.}
		\label{fig:delta_Rstr1}
	\end{figure}

\begin{figure}[H]
		\centering
		\begin{subfigure}[b]{0.38\linewidth}
			\includegraphics[width=\linewidth]{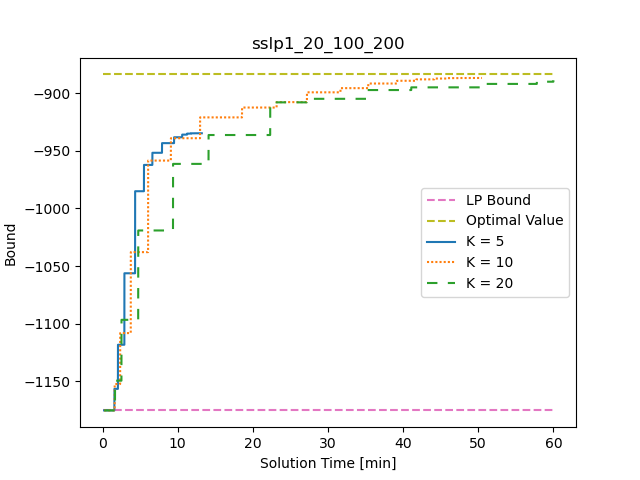}
		\end{subfigure}
		\begin{subfigure}[b]{0.38\linewidth}
			\includegraphics[width=\linewidth]{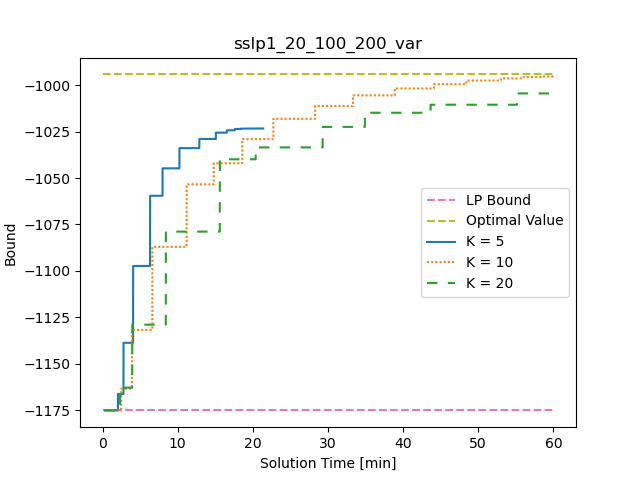}
		\end{subfigure}
		\begin{subfigure}[b]{0.38\linewidth}
			\includegraphics[width=\linewidth]{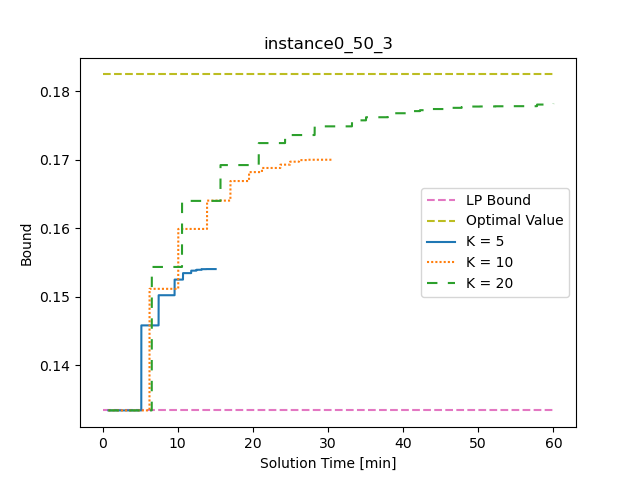}
		\end{subfigure}
		\caption{Convergence profile of {\rstrone} on an SSLP instance (top left), {an SSLPV
				instance (top right),} and a SNIP instance (bottom) with $\delta=50\%$ and varying $K$ values.}
		\label{fig:Rstr1_K}
	\end{figure}

\begin{figure}[H]
		\centering
		\begin{subfigure}[b]{0.38\linewidth}
			\includegraphics[width=\linewidth]{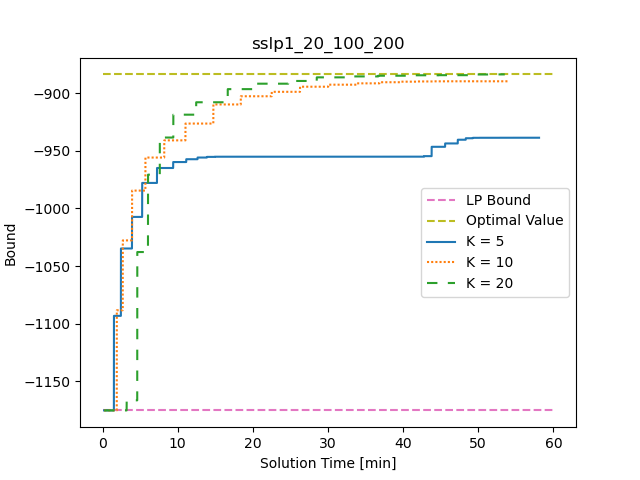}
		\end{subfigure}
		\begin{subfigure}[b]{0.38\linewidth}
			\includegraphics[width=\linewidth]{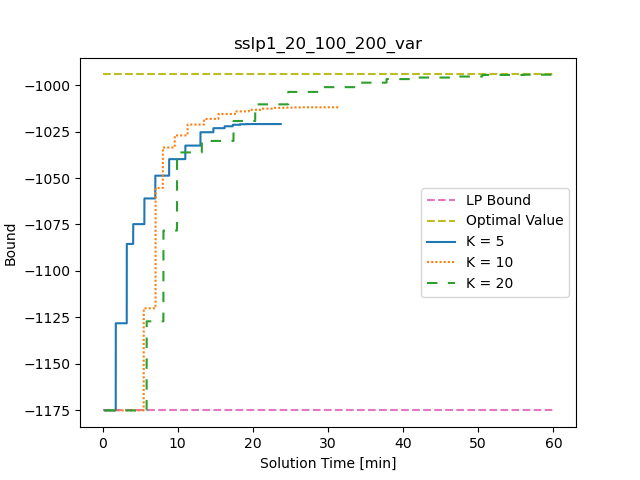}
		\end{subfigure}
		\begin{subfigure}[b]{0.38\linewidth}
			\includegraphics[width=\linewidth]{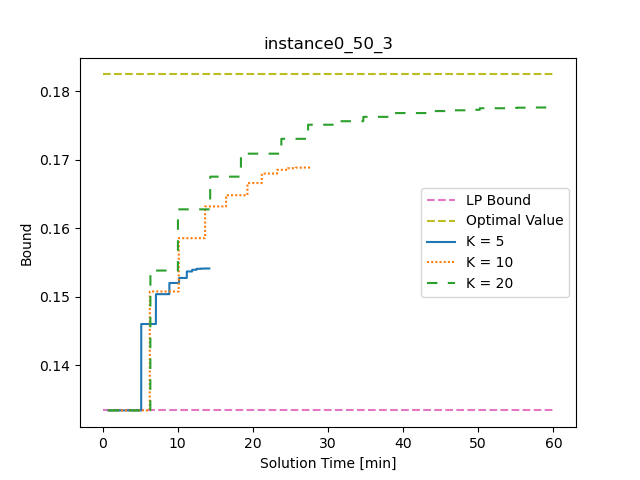}
		\end{subfigure}
		\caption{Convergence profile of {\rstrtwo} on an SSLP instance (top left), {an SSLPV
				instance (top right),} and a SNIP instance (bottom) with $\delta=50\%$ and varying $K$ values.}
		\label{fig:Rstr2_K}
	\end{figure}

\begin{figure}[H]
		\centering
		\begin{subfigure}[b]{0.38\linewidth}
			\includegraphics[width=\linewidth]{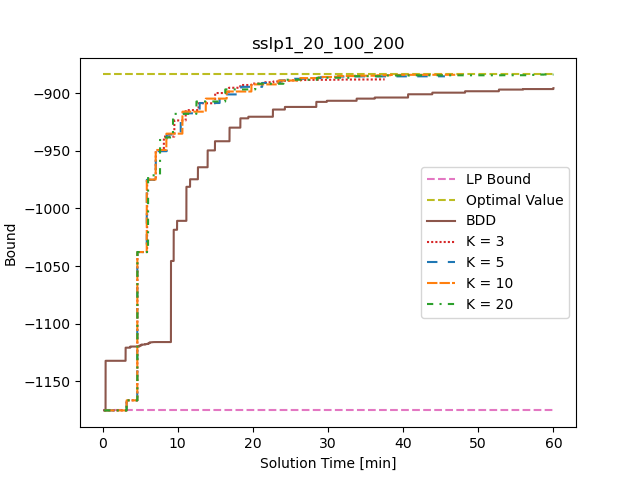}
		\end{subfigure}
		\begin{subfigure}[b]{0.38\linewidth}
			\includegraphics[width=\linewidth]{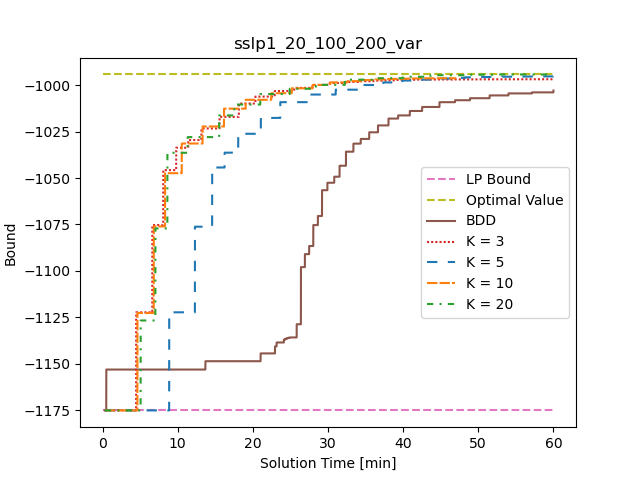}
		\end{subfigure}
		\begin{subfigure}[b]{0.38\linewidth}
			\includegraphics[width=\linewidth]{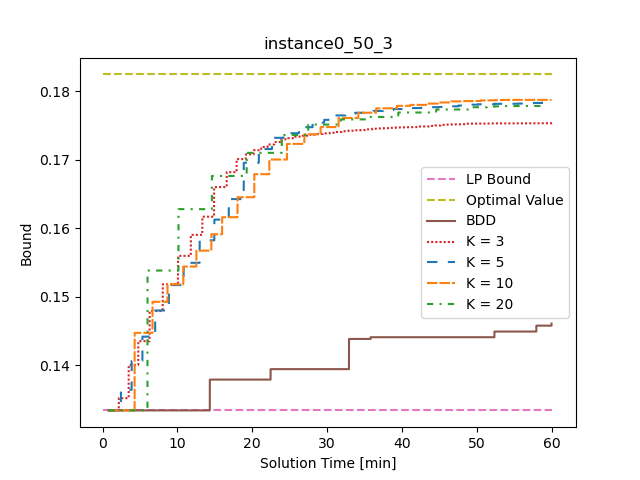}
		\end{subfigure}
		\caption{Convergence profile of {\rstrmip} on an SSLP instance (top left), {an SSLPV
				instance (top right),} and an SNIP instance (bottom) with $\delta=50\%$ and varying $K$ values and comparison with {\BDD}.}
		\label{fig:RstrMIP_K}
	\end{figure}

	Figures \ref{fig:delta_Rstr1} - \ref{fig:RstrMIP_K} present the evolution of the progress in the
	lower bound over time on one representative example SSLP instance and one representative example SNIP instance using
	various cut generation approaches and parameters.
	
	\begin{figure}[hbt!]
		\centering
		\begin{subfigure}[b]{0.6\linewidth}
			\includegraphics[width=\linewidth]{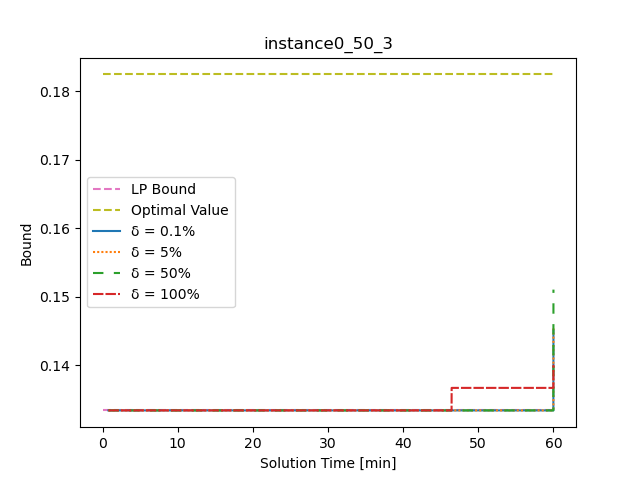}
		\end{subfigure}
		\caption{Convergence profile of an SNIP instance (instance0\_50\_3) with varying $\delta$ values.}
		\label{fig:delta_Exact_SNIP}
	\end{figure}
	
	
	From Figure \ref{fig:RstrMIP_K} we see that for the example SSLP instance, {\BDD} does not make much progress in the
	first 10 minutes, but then makes steady improvement. The shift occurs when {\BDD} switches from its initial phase of
	generating only strengthened Benders cuts to the second phase of heuristically generating Lagrangain cuts. In the
	initial phase 
	{\BDD} the strengthened Benders cuts are not strong enough to substantially improve the
	bound after the first iteration, so that significant bound progress only occurs again after switching to the second
	phase. We experimented with making this transition quicker (or even skipping it altogether) but we found that this led to poor
	performance of {\BDD}. It  appears that the first phase is important for building an inner approximation that is used for
	separation of Lagrangian cuts in the next phase. The convergence profiles of {\BDD} on the SNIP instance and the SSLPV instance clearly show its
	slow convergence relative to all variations of our proposed restricted Lagrangian cut generation approach.
	
	\bibliographystyle{ieeetr}
	\bibliography{joc/ref}
\end{document}

%% file: joc/main.tex
\section{Introduction}\label{sec:intro}
We study methods for solving two-stage stochastic integer programs (SIPs) with general mixed-integer first-stage and
second-stage variables. Two-stage stochastic programs are used to model problems with uncertain data, where a decision
maker first decides the values of first-stage decision variables, then observes the values of the uncertain data, and
finally decides the values of second-stage decision variables. The objective is to minimize the sum of the first-stage
cost and the expected value of the second-stage cost, where the expected value is taken with respect to the distribution
of the uncertain data. Each realization of the uncertain data is called a scenario. Assuming the uncertainty is modeled
with a finite set of scenarios $S$, a two-stage SIP can be formulated as follows:\begin{equation}\label{SIP}
z_{\text{IP}}=\min_x\Big\{c^\top x+\sum_{s\in S}p_sQ_s(x):Ax\geq b,x\in X \Big\},
\end{equation}
where $c\in\bR^n$, and for each $s\in S$, $p_s$ denotes the probability of scenario $s$ and $Q_s:\bR^n\rightarrow\bR\cup\{+\infty\}$ is the recourse function of scenario $s$ defined as\begin{equation}\label{def:Qs}
Q_s(x)=\min_y\{(q^s)^\top y:\ W^sy\geq h^s-T^sx,y\in Y\}.
\end{equation}
Here $X\subseteq\bR^n$ and $Y\subseteq\bR^{n_y}$ denote the integrality restrictions on some or all variables of $x$ and
$y$, respectively. Sign restrictions and variable bounds, if present, are assumed to be included in the constraints $Ax \geq b$ or $W^s y
\geq h^s - T^s x$. Matrices $W^s,T^s$ and vectors $q^s,h^s$ are the realization of the uncertain data associated with
scenario $s\in S$. We do not assume relatively complete recourse, i.e., it is possible that the recourse problem
\eqref{def:Qs} is infeasible for some $s \in S$ and $x \in X$  satisfying $Ax\geq b$, in which case 
$Q_s(x)=+\infty$ by convention.

Two-stage SIPs can also be written in the extensive form:\begin{equation}
\begin{aligned}\label{ext_form}
\min_{x,y^s}\ & c^\top x+\sum_{s\in S}p_s(q^s)^\top y^s\\
\text{s.t. }\ &Ax\geq b,x\in X,\\
&W^sy^s\geq h^s-T^sx^s,y^s\in Y,s\in S.
\end{aligned}
\end{equation}
From this perspective, SIPs are essentially large-scale mixed-integer programs (MIPs) with special block structure.
Directly solving \eqref{ext_form} as a MIP can be difficult when $|S|$ is large. Therefore, significant research has
been devoted to the study of decomposition methods for solving SIP problems. Benders decomposition
\pv{\citep{bnnobrs1962partitioning,van1969shaped}}\gv{\cite{bnnobrs1962partitioning,van1969shaped}} and dual
decomposition \pv{\citep{caroe1999dual}}\gv{\cite{caroe1999dual}} are the two most commonly used decomposition methods
for solving SIPs. In Benders decomposition, the second-stage variables are projected out from the Benders master model
and linear programming (LP) relaxations of scenario subproblems are solved to add cuts in the master model. In dual
decomposition, a copy of the first-stage variables is created for each scenario and constraints are added to require the
copies to be equal to each other. A Lagrangian relaxation is then formed by relaxing these so-called nonanticipativity constraints. Solving the
corresponding Lagrangian dual problem requires solving single-scenario MIPs to provide function values and supergradients
of the Lagrangian relaxation. A more detailed description of dual decomposition for SIPs is given
in Section \ref{sect:dd}. Usually dual decomposition generates a much stronger bound than Benders decomposition but the
bound takes much longer to compute. \gv{Rahmaniani et al.
\cite{rahmaniani2020benders}}\pv{\cite{rahmaniani2020benders}} proposed Benders dual decomposition (BDD) in which 
Lagrangian cuts generated by solving single-scenario MIPs identical to the subproblems in dual decomposition are added
to the Benders formulation to strengthen the relaxation.
Similarly, \gv{Li and Grossmann
\cite{li2018improved}}\pv{\cite{li2018improved}} develop a Benders-like decomposition algorithm implementing both
Benders cuts and Lagrangian cuts for solving convex mixed 0-1 nonlinear stochastic programs. 
While representing an interesting hybrid between the Benders
and dual decomposition approaches, the time spent generating the Lagrangian cuts remains a limitation in these
approaches. In this paper, we extend this line of work by investigating methods to generate Lagrangian cuts more efficiently.

This work more broadly contributes to the significant research investigating the use of cutting planes to strengthen the Benders
model. \gv{Laporte and Louveaux \cite{laporte1993integer}}\pv{\cite{laporte1993integer}} propose integer L-shaped cuts
for SIPs with pure binary first-stage variables. \gv{Sen and Higle \cite{sen2005c}}\pv{\cite{sen2005c}} and \gv{Sen and
Sherali \cite{sen2006decomposition}}\pv{\cite{sen2006decomposition}} apply disjunctive programming techniques to
convexify the second-stage problems for SIPs with pure binary first-stage variables. \gv{Ntaimo
\cite{ntaimo2013fenchel}}\pv{\cite{ntaimo2013fenchel}} investigates the use of Fenchel cuts added to strengthen the
subproblem relaxations. \gv{Gade et al. \cite{gade2014decomposition} and Zhang and K{\"u}{\c{c}}{\"u}kyavuz
\cite{zhang2014finitely}}\pv{\cite{gade2014decomposition} and \cite{zhang2014finitely}} demonstrate how Gomory cuts can
be used within a Benders decomposition method for solving SIPs with pure integer first-stage and second-stage variables, leading to
a finitely convergent algorithm.
\gv{Qi and Sen \cite{qi2017ancestral}}{\pv{\cite{qi2017ancestral}} use cuts valid for multi-term disjunctions
introduced in 
\cite{chen2011finite} to derive an algorithm for SIPs with mixed-integer recourse.
\gv{{van} der Laan and Romeijnders \cite{van2020converging}}{\pv{\cite{van2020converging}} propose a new class of cuts,
scaled cuts, which integrates previously generated Lagrangian cuts into a cut-generation model, again leading to a finitely
convergent cutting-plane method for SIPs with mixed-integer recourse. 
\gv{Bodur et al.
\cite{bodur2017strengthened}}\pv{\cite{bodur2017strengthened}} compare the strength of split cuts derived in the Benders
master problem space to those added in the scenario LP relaxation space. 
\gv{Zou et al.
\cite{zou2019stochastic}}\pv{\cite{zou2019stochastic}} propose strengthened Benders cuts and Lagrangian cuts to solve
pure binary multistage SIPs. 

There are also other approaches for solving SIPs that do not rely on a Benders model. \gv{Lulli and Sen \cite{lulli2004branch}}\pv{\cite{lulli2004branch}} develop a
column generation based algorithm for solving multistage SIPs. \gv{Lubin et al.
\cite{lubin2013parallelizing}}\pv{\cite{lubin2013parallelizing}} apply proximal bundle methods to parallelize the dual
decomposition algorithm. \gv{Boland et al.  \cite{boland2018combining}}\pv{\cite{boland2018combining}} and \gv{Guo et
al. \cite{phbounds:orl15}}\pv{\cite{phbounds:orl15}} investigate the use of progressive hedging to calculate the
Lagrangian dual bound. \gv{Kim and Dandurand \cite{dandurandscalable}}\pv{\cite{dandurandscalable}} propose a new
branching method for dual decomposition of SIPs. \gv{Kim et al.
\cite{kim2019asynchronous}}\pv{\cite{kim2019asynchronous}} apply asynchronous trust-region methods to solve the
Lagrangian dual used in dual decomposition.

The goal of this paper is to develop an effective way of generating Lagrangian cuts to add to a Benders model. The main contributions of our work are summarized as follows.\begin{enumerate}
	\item We propose a normalization for generating Lagrangian cuts similar to the one used \pv{by
	\cite{fischetti2010note}}\gv{in \cite{fischetti2010note}} for separating Benders cuts. This normalization can be used to construct a Lagrangian cut separation problem different from the one used \pv{by \cite{zou2019stochastic} and \cite{rahmaniani2020benders}}\gv{in \cite{zou2019stochastic} and \cite{rahmaniani2020benders}}.
	\item We propose methods for accelerating the generation of Lagrangian cuts, including solving the cut generation
	problem in a restricted subspace, and using a MIP approximation to identify a promising restricted subspace. Numerical
	results indicate that these approaches lead to significantly faster bound improvement from Lagrangian cuts.
	\item We conduct an extensive numerical study on three classes of two-stage SIPs. We compare the impact of different strategies for generating Lagrangian cuts. Computational results are also given for using our method within branch-and-cut as an exact solution method. We find that this method has potential to outperform both a pure Benders approach and a pure dual decomposition approach.
\end{enumerate}
Finally, we note that while we focus on new techniques for generating Lagrangian cuts for use in solving two-stage SIPs,
these techniques can be directly applied to multi-stage SIPs by integrating them within the stochastic dual dynamic
integer programming (SDDIP) algorithm \pv{by \cite{zou2019stochastic}}\gv{in \cite{zou2019stochastic}}.  

\section{Preliminaries}	\label{sec:Pre}
\subsection{Branch-and-Cut Based Methods}\label{sect:bnc}
Problem \eqref{SIP} can be reformulated as the following Benders model:\begin{equation}\label{Benders}
\min_{x,\theta_s}\Big\{c^\top x+\sum_{s\in S}p_s\theta_s:(x,\theta_s)\in E^s, s\in S \Big\},
\end{equation}
where for each $s\in S$, $E^s$ contains the first-stage constraints and the epigraph of $Q_s$, i.e.,\begin{equation}\label{def:Es}
E^s=\{(x,\theta_s)\in X\times\bR:Ax\geq b,\theta_s\geq Q_s(x) \}.
\end{equation}
In a typical approach to solve \eqref{Benders}, each recourse function $Q_s$ is replaced by a cutting-plane underestimate $\hat{Q}_s$ which
creates a relaxation and is dynamically updated. In each iteration, an approximate problem defined by the current
underestimate is solved to obtain a candidate solution and then a cut generation problem is solved to update the cutting
plane approximation if necessary. This process is repeated until no cuts are identified.

We review two types of valid inequalities for $E^s$. The first collection of cuts are Benders cuts
\pv{\citep{bnnobrs1962partitioning,van1969shaped}}\gv{\cite{bnnobrs1962partitioning,van1969shaped}}. Benders cuts are
generated based on the LP relaxation of the problem \eqref{def:Qs} defining $Q_s(x)$. Given a candidate solution
$\hat{x}$, the LP relaxation of the recourse problem is solved:\begin{equation}\label{BendersSub}
\min_{y}\{(q^s)^\top y:W^sy\geq h^s-T^s\hat{x}\}.
\end{equation}
Let $\mu$ be an optimal dual solution of problem \eqref{BendersSub}. Based on LP duality, the following {\it Benders
cut} is valid for $E^s$:\begin{equation}\label{BenCut}
\mu^\top T^sx+\theta_s\geq \mu^\top h^s.
\end{equation}
If the SIP has continuous recourse, i.e., $Y=\bR^{n_y}$, then \eqref{BenCut} is tight at $\hat{x}$, i.e., $Q_s(\hat{x})= -\mu^\top T^s\hat{x}+\mu^\top h^s.$ The cutting-plane model $\hat{Q}_s$ is often constructed by iteratively adding Benders cuts until the lower bound converges to the LP relaxation bound. Benders cuts are sufficient to provide convergence for solving SIPs with continuous recourse \pv{\citep{kall1994stochastic}}\gv{\cite{kall1994stochastic}}.

Another useful family of cuts is the integer L-shaped cuts introduced \pv{by \cite{laporte1993integer}}\gv{in \cite{laporte1993integer}}. These cuts are valid
only when the first-stage variables are binary, i.e., $X=\{0,1\}^n$, but can be applied even when the second-stage
includes integer decision variables. In this case, given $\hat{x}\in\{0,1\}^n$ and a value $L_s$ such that $Q_s(x)\geq
L_s$ for all feasible $x$, the following {\it integer L-shaped cut} is a valid inequality for $E^s$:\begin{equation}\label{IntLShaped}
\theta_s\geq Q_s(\hat{x})-(Q_s(\hat{x})-L_s)\Big(\sum_{i:\hat{x}_i=1}(1-x_i)+\sum_{i:\hat{x}_i=0}x_i\Big).
\end{equation}

A standard branch-and-cut algorithm using Benders cuts and integer L-shaped cuts is described in \pv{Section
S.\ref{BnC} of the online supplement}\gv{Appendix}.
The branch-and-cut algorithm can be implemented using a lazy constraint callback in modern MIP solvers, which allows the
addition of Benders or integer L-shaped cuts when the solver encounters a solution $(\hat{\theta},\hat{x})$ with
$\hat{x} \in X$ (i.e., $\hat{x}$ satisfies any integrality constraints) but for which $\hat{\theta}_s < Q_s(\hat{x})$
for some $s \in S$.  These two classes of cuts are sufficient to guarantee convergence for SIPs with continuous recourse
or pure binary first-stage variables. However, the efficiency of the algorithm depends significantly on the strength of
the cutting-plane models $\hat{Q}_s$'s. Given poor relaxations of the recourse functions, the branch-and-bound search
may end up exploring a huge number of nodes, resulting in a long solution time. As discussed in Section \ref{sec:intro}
many methods have been proposed to strengthen the model $\hat{Q}_s$ in order to accelerate the algorithm.

We note that the validity of Lagrangian cuts does not require either continuous recourse or pure binary first-stage
variables. However, outside of those settings, additional cuts or a specialized branching scheme would be required to
obtain a convergent algorithm. We refer the reader to, e.g.,
\cite{ahmed2004finite,qi2017ancestral,van2020converging,zhang2014finitely} for examples of methods that
can be used to obtain a finitely convergent algorithm in other settings. Lagrangian cuts could potentially be added to
enhance any of these approaches.

\subsection{Dual Decomposition}\label{sect:dd}
In dual decomposition \pv{\citep{caroe1999dual}}\gv{\cite{caroe1999dual}}, a copy $x^s$ of the first-stage variables $x$ is created for each scenario
$s \in S$ which yields a reformulation of \eqref{SIP}:\begin{align*}
\min_{x,x^s,y^s} &\sum_{s\in S}p_s(c^\top x^s+(q^s)^\top y^s)\\
&\begin{array}{ll}
\hspace*{-0.35in}\text{s.t. }\:\:\: Ax^s\geq b,&s\in S,\\
T^sx^s+W^sy^s\geq h^s, &s\in S,\\
x^s\in X,y^s\in Y, &s\in S,\\
x^s=x, &s\in S.
\end{array}
\end{align*}
Lagrangian relaxation is applied to the {\it nonanticipativity constraints} $x^s=x$ in this formulation with multipliers $\lambda^s$ for $s\in S$, which gives the following Lagrangian relaxation problem:\begin{equation}\label{DD}
\begin{aligned}
z(\lambda)=\min_{x,x^s,y^s} &\sum_{s\in S}p_s(c^\top x^s+(q^s)^\top y^s)+\sum_{s\in S}p_s(\lambda^s)^\top (x^s-x),\\
\text{s.t. } &(x^s,y^s)\in K^s,\quad s\in S,
\end{aligned}
\end{equation}
where $K^s:=\{x \in X, y \in Y : Ax\geq b, T^sx+W^sy\geq h^s\}$ for each $s\in S$. Throughout this paper we assume that $K^s$ is nonempty for each $s\in S$ 

The nonanticipative Lagrangian dual problem is defined as $z_D=\max_{\lambda} z(\lambda)$.
Note that constraint $\sum_{s\in S}p_s\lambda^s=0$ is implicitly enforced in $\max_\lambda z(\lambda)$ since
$z(\lambda)=-\infty$ when $\sum_{s\in S}p_s\lambda^s\neq0$. Under this condition, the variables $x$ in \eqref{DD} can be
eliminated, and hence \eqref{DD} can be solved by solving a separate optimization for each scenario $s \in S$. The nonanticipative Lagrangian dual problem \begin{displaymath}
z_D=\max_{\lambda}\Big\{z(\lambda):\sum_{s\in S}p_s\lambda^s=0\Big\}
\end{displaymath}
is a convex program which gives a lower bound on the optimal value of \eqref{SIP}. The following theorem provides a primal characterization of the bound $z_D$.
\begin{thm} [\pv{\cite{caroe1999dual}}\gv{Car{\o}e and Schultz \cite{caroe1999dual}}]\label{thm:LD_bound}
	The following equality holds:\begin{displaymath}
	z_D=\min_{x,y^s}\Big\{c^\top x+\sum_{s\in S}p_s(q^s)^\top y^s:(x,y^s)\in \conv(K^s),s\in S \Big\}.
	\end{displaymath}
\end{thm}

Empirical evidence \pv{\citep[e.g.][]{schutz2009supply,solak2010optimization}}\gv{(e.g.,
\cite{schutz2009supply,solak2010optimization})} shows that the Lagrangian dual bound $z_D$ is often a tight lower bound
on the optimal objective value of \eqref{def:Qs}. Such a bound can be used in a branch and bound algorithm and to
generate feasible solutions using heuristics \pv{\citep[e.g.][]{caroe1999dual,dandurandscalable}}\gv{(see, e.g.,
\cite{caroe1999dual,dandurandscalable})}. However, the Lagrangian dual problem can be hard to solve due to its large
size (with $n|S|$ variables) and the difficulty of solving MIP subproblems to evaluate $z(\lambda)$ and obtain supergradients of $z(\cdot)$ at a point $\lambda$.

\subsection{Benders Dual Decomposition}
The BDD algorithm \pv{\citep{rahmaniani2020benders}}\gv{\cite{rahmaniani2020benders}} is a version of Benders decomposition that uses strengthened Benders cuts and Lagrangian cuts to tighten the Benders model. In particular, they solve two types of scenario MIPs to generate optimality and feasibility cuts for $E^s$:\begin{enumerate}
	\item Given any $\lambda\in \bR^n$, let $\big(\xslambda,\yslambda\big)$ be an optimal solution of\begin{displaymath}
	\min_{x,y}\{\lambda^\top x+(q^s)^\top y:(x,y)\in K^s \}.
	\end{displaymath}
	Then the following {\it Lagrangian optimality cut} is valid for $E^s$:\begin{equation}\label{ineq:Lag}
	\lambda^\top (x-\xslambda)+\theta_s\geq (q^s)^\top \yslambda.
	\end{equation}
	\item Given any $\lambda\in \bR^n$, let $\big(\hat{x}^s_\lambda,\hat{y}^s_\lambda,\hat{u}^s_\lambda,\hat{v}^s_\lambda\big)$ be an optimal solution of\begin{displaymath}
	\min_{x,y,u,v}\{\mathbbm{1}^\top v+\mathbbm{1}^\top u+\lambda^\top x:Ax+u\geq b,T^sx+W^sy+v\geq h^s,x\in X,y\in Y \}.
	\end{displaymath}
	Then the following {\it Lagrangian feasibility cut} is valid for $E^s$:\begin{equation}\label{ineq:fea}
	\lambda^\top (x-\hat{x}^s_\lambda) \geq\mathbbm{1}^\top \hat{v}^s_\lambda+\mathbbm{1}^\top \hat{u}^s_\lambda.
	\end{equation}
\end{enumerate}
The BDD algorithm generates both types of Lagrangian cuts by heuristically solving a Lagrangian cut generation problem.
Numerical results from \cite{rahmaniani2020benders} show that Lagrangian cuts are able to close significant gap at the
root node for a variety of SIP problems. Our goal in this work is to provide new methods for quickly finding strong Lagrangian cuts.

\section{Selection and Normalization of Lagrangian Cuts}

We begin by introducing an alternative view of the Lagrangian cuts \eqref{ineq:Lag} and \eqref{ineq:fea} which leads to
a different cut generation formulation. Let $\bQs:\bR^n\times\bR_+ \rightarrow \bR$ be defined as 
\begin{align}\label{subIP}
\bQs(\pi,\pi_0)=\min_x&\bigl\{\pi^\top x+\pi_0Q_s(x): Ax \geq b, x \in X \bigr\}=\min_{x,y}\bigl\{\pi^\top x+\pi_0(q^s)^\top y:(x,y)\in K^s\bigr\}\end{align}
for $(\pi,\pi_0) \in\bR^n\times\bR_+$. 
According to \eqref{subIP}, for any $(\pi,\pi_0)\in\bR^n\times\bR_+$, the following inequality is valid for set $E^s$ defined in \eqref{def:Es}:\begin{equation}\label{cut}
\pi^\top x+\pi_0\theta_s\geq \bQs(\pi,\pi_0).
\end{equation}

We next demonstrate that the set of all cuts of the form \eqref{cut} is equivalent to set of Lagrangian cuts \eqref{ineq:Lag} and
\eqref{ineq:fea}, and therefore we refer to cuts of the form \eqref{cut} as Lagrangian cuts. 
For any $\Pi_s\subseteq \bR^n\times\bR_+$, let $P_s(\Pi_s)$ denote the set defined by all cuts \eqref{cut} with coefficients restricted on $\Pi_s$, i.e.,\begin{displaymath}
\PsPis:=\{(x,\theta_s): \pi^\top x+ \pi_0\theta_s\geq \bQs(\pi,\pi_0)\text{ for all }(\pi,\pi_0)\in\Pi_s \}.
\end{displaymath}

\begin{prop}\label{prop:strength}
	Let $\Pi_s\subseteq\bR^{n+1}$ be a neighborhood of the origin and let $\Pi_s^*=\Pi_s\cap(\bR^n\times\bR_+)$. Define $\bar{P}_s:=\{(x,\theta_s):\text{\eqref{ineq:Lag} and \eqref{ineq:fea} hold for all }\lambda\in\bR^n\}$. Then $\PsPistar=\bar{P}_s$.
\end{prop}
\pv{\proof{Proof.}}
\gv{\begin{proof}}
We first show that $\PsPistar\subseteq \bar{P}_s$. Let $(x,\theta_s)\in \PsPistar$. We show that $(x,\theta_s)$
satisfies \eqref{ineq:Lag} and \eqref{ineq:fea} for all $\lambda\in\bR^n$. Let $\lambda\in\bR^n$, and choose
$\rho^1_\lambda$ and $\rho^2_\lambda$ large enough such that
$(\lambda/\rho^1_\lambda,1/\rho^1_\lambda),(\lambda/\rho^2_\lambda,0)\in\Pi_s^*$, and thus
$(x,\theta_s)$ satisfies \eqref{cut} using these values for $(\pi,\pi_0)$.
Multiplying inequality \eqref{cut} with $(\pi,\pi_0)=(\lambda/\rho^1_\lambda,1/\rho^1_\lambda)\in\Pi_s^*$ by $\rho_\lambda^1$ gives\begin{displaymath}
\lambda^\top x+\theta_s \geq \rho^1_\lambda\bQs(\lambda/\rho^1_\lambda,1/\rho^1_\lambda)=\bQs(\lambda,1)=\lambda^\top \xslambda+(q^s)^\top \yslambda,
\end{displaymath}
where $(\xslambda,\yslambda)$ is as defined in \eqref{ineq:Lag} and 
where we have used the observation that the function $\bQs$ is positively homogeneous.
Subtracting $\lambda^\top \xslambda$ from both sides yields that $(x,\theta_s)$ satisfies \eqref{ineq:Lag}. On the other hand, inequality \eqref{cut} with $(\pi,\pi_0)=(\lambda/\rho^2_\lambda,0)\in\Pi_s^*$ implies
\begin{align*}
(\lambda/\rho^2_\lambda)^\top x &\geq\bQs(\lambda/\rho^2_\lambda,0) \\
&= \min_{x,y}\{\lambda^\top x:Ax\geq b,T^sx+W^sy\geq h,x\in X,y\in Y \}/\rho^2_\lambda\\
&\geq\min_{x,y,u,v}\{\mathbbm{1}^\top v+\mathbbm{1}^\top u+\lambda^\top x:Ax+u\geq b,T^sx+W^sy+v\geq h^s,x\in X,y\in Y \}/\rho^2_\lambda\\
&=(\lambda^\top \hat{x}^s_\lambda+\mathbbm{1}^\top \hat{v}^s_\lambda+\mathbbm{1}^\top \hat{u}^s_\lambda)/\rho^2_\lambda.
\end{align*}
Multiplying both sides by $\rho_\lambda^2$ and subtracting $\lambda^\top \hat{x}^s_\lambda$ from both sides gives inequality \eqref{ineq:fea}. Since $\lambda$ can be arbitrary, it implies $\PsPistar\subseteq \bar{P}_s$.

Next let $(x,\theta_s)\in\bar{P}_s$. We aim to show $(x,\theta_s)\in\PsPistar$. Let $(\pi,\pi_0)\in\Pi_s^*$. If
$\pi_0>0$, then \eqref{cut} is satisfied by scaling \eqref{ineq:Lag} with $\lambda=\pi/\pi_0$. We can similarly show that
$\pi^\top x+\pi_0'\theta_s\geq\bQs(\pi,\pi_0')$ for any $\pi_0' > 0$. Now assume $\pi_0=0$. By
Theorem 7.1 of \cite{rockafellar1970convex}, $\bQs$ is lower-semicontinuous, and hence \begin{displaymath}
\pi^\top  x =\pi^\top  x + \liminf_{\pi_0'\rightarrow 0}\pi_0'\theta_s\geq \liminf_{\pi_0'\rightarrow 0}\bQs(\pi,\pi_0')\geq\bQs(\pi,0),
\end{displaymath}
i.e., \eqref{cut} is satisfied. Therefore, we have $\bar{P}_s\subseteq \PsPistar$.
\pv{\Halmos}
\pv{\endproof}
\gv{\end{proof}}

Since Proposition \ref{prop:strength} holds for any neighborhood $\Pi_s$ of the origin, it allows flexibility to choose
a normalization on the Lagrangian cut coefficients when generating cuts. We discuss the choice of such a normalization in Section \ref{sect:alg}.

Let $z_{\text{LC}}$ denote the bound we obtain after adding all Lagrangian cuts, i.e.,\begin{displaymath}
z_{\text{LC}}:=\min_{x,\theta_s}\Bigl\{c^\top x+\sum_{s\in S}p_s\theta_s: (x,\theta_s)\in\Ps{\bR^n\times\bR_+} ,s\in S \Bigr\}.
\end{displaymath} 
Since $\bQs$ is positively homogeneous and the set defined by an inequality is invariant to any positive scaling of the
inequality, we also have
\begin{displaymath}
z_{\text{LC}}=\min_{x,\theta_s}\Bigl\{c^\top x+\sum_{s\in S}p_s\theta_s: (x,\theta_s)\in \PsPistar,s\in S \Bigr\}
\end{displaymath} 
for any $\Pi_s^*=\Pi_s\cap(\bR^n\times\bR_+)$ such that $\Pi_s$ is a neighborhood of the origin. \gv{Rahmaniani et al.
\cite{rahmaniani2020benders}}\pv{\cite{rahmaniani2020benders}} showed that for $|S|=1$, adding all Lagrangian cuts of
the form \eqref{ineq:Lag} and feasibility cuts \eqref{ineq:fea} yields a Benders master problem with optimal objective value equal to $z_{\text{IP}}$.
This is not true in general when $|S| > 1$. 
However, \gv{Li and Grossmann \cite{li2018improved}}\pv{\cite{li2018improved}} showed that
$z_{\text{LC}}\geq z_D$. We show in the following theorem that this is an equality, i.e., $z_{\text{LC}}=z_D$.
Since $z_{\text{IP}}=z_D$ according to Theorem \ref{thm:LD_bound} when $|S|=1$, this result generalizes the result \pv{by \cite{rahmaniani2020benders}}\gv{in \cite{rahmaniani2020benders}}.  
\begin{thm}\label{thm:eqv}
	The equality $z_{\text{LC}}=z_D$ holds.
\end{thm}
\pv{\proof{Proof.}}
\gv{\begin{proof}}
Although $z_{\text{LC}} \geq z_D$ has been shown in \cite{li2018improved}, we show both directions for completeness.

	By Theorem \ref{thm:LD_bound}, \begin{multline*}
	z_D=\min_{x,y^s}\Big\{c^\top x+\sum_{s\in S}p_s(q^s)^\top y^s:(x,y^s)\in \conv(K^s),s\in S \Big\}\\
	=\min_{x,\theta_s}\Big\{c^\top x+\sum_{s\in S}p_s\theta_s:(x,\theta_s)\in U^s,s\in S\Big\},
	\end{multline*}
	where $U^s:=\proj_{(x,\theta_s)}\{(x,y^s,\theta_s):\theta_s\geq (q^s)^\top y^s, (x,y^s)\in \conv(K^s)\}$. Therefore, we only need to show that $U^s=\Ps{\bR^n\times\bR_+}$ for each $s\in S$.
	
	On the one hand, let $s\in S$ and $(\bar{x},\bar{\theta}_s)\in U^s$. Then there exists $\bar{y}^s$ such that $\bar{\theta}_s\geq (q^s)^\top \bar{y}^s$ and $(\bar{x},\bar{y}^s)\in\conv(K^s)$. Then by definition of $\bQs$, we have\begin{multline*}
	\pi^\top \bar{x}+\pi_0\bar{\theta}_s\geq\pi^\top \bar{x}+\pi_0(q^s)^\top \bar{y}^s
	\geq\min_{x,y}\{\pi^\top x+\pi_0(q^s)^\top y:(x,y)\in\conv(K^s) \}=\bQs(\pi,\pi_0)
	\end{multline*}
	for all $(\pi,\pi_0)\in\bR^n\times\bR_+$, which implies that $(\bar{x},\bar{\theta}_s)\in \Ps{\bR^n\times\bR_+}$.
	
	On the other hand, let $s\in S$ and $(\bar{x},\bar{\theta}_s)\notin U^s$. Since $U^s$ is a polyhedron, which is convex and closed, there exists a hyperplane strictly separating $(\bar{x},\bar{\theta}_s)$ from $U^s$. In other words, there exists $u\in\bR^n$ and $v\in \bR$ such that $u^\top \bar{x}+v\bar{\theta}_s<\min_{x,\theta_s}\{u^\top x+v\theta_s:(x,\theta_s)\in U^s \}$. Note that we have $v\geq 0$ in this case since $\min_{x,\theta_s}\{u^\top x+v\theta_s:(x,\theta_s)\in U^s \}=-\infty$ if $v<0$. Then\begin{multline*}
	u^\top \bar{x}+v\bar{\theta}_s<\min_{x,\theta_s}\{u^\top x+v\theta_s:(x,\theta_s)\in U^s \}\\
	=\min_{x,y}\{u^\top x+v(q^s)^\top y:(x,y)\in \conv(K^s) \}=\bQs(u,v),
	\end{multline*}
	i.e., $(\bar{x},\bar{\theta}_s)$ violates the Lagrangian cut $u^\top \bar{x}+v\bar{\theta}_s\geq\bQs(u,v)$. Therefore, $(\bar{x},\bar{\theta}_s)\notin \Ps{\bR^n\times\bR_+}$.
\pv{\Halmos}
\pv{\endproof}
\gv{\end{proof}}

We next consider the problem of separating a point from $\PsPistar$ using Lagrangian cuts. Given a candidate solution $(\hat{x},\hat{\theta}_s)$ and a convex and compact $\Pi_s^*$, the problem of separating $(\hat{x},\hat{\theta}_s)$ from $\PsPistar$ can be formulated as the following bounded convex optimization problem\begin{equation}
\label{separation}
\max_{\pi,\pi_0}\{\bQs(\pi,\pi_0)-\pi^\top \hat{x}-\pi_0^\top \hat{\theta}_s:(\pi,\pi_0)\in\Pi_s^* \}.
\end{equation}
Problem \eqref{separation} can be solved by bundle methods given an oracle for evaluating the function value and a
supergradient of $Q_s^*$. The function value and a supergradient of $\bQs$ can be evaluated at any
$(\pi,\pi_0)\in\bR^n\times\bR_+$ by solving the single-scenario MIP \eqref{subIP}.
Let $(x^*,y^*)$ be an optimal solution of \eqref{subIP}. Then we have
$\bQs(\pi,\pi_0)=\pi^\top x^*+\pi_0(q^s)^\top y^*$ and $(x^*,(q^s)^\top y^*)$ is a supergradient of $\bQs$ at $(\pi,\pi_0)$.

Although \eqref{separation} is a convex optimization problem, the requirement to solve a potentially large number of MIP problems of the form
\eqref{subIP} to evaluate function values and supergradients of $\bQs$  to find a single Lagrangian cut can make the use of such cuts computationally prohibitive. Numerical results \pv{by}\gv{in}
\cite{zou2019stochastic} indicate that although Lagrangian cuts can significantly improve the LP relaxation bound for
many SIP problems, the use of exact separation of Lagrangian cuts does not improve overall solution time. In an attempt
to reduce the time required to separate Lagrangian cuts, \pv{\cite{rahmaniani2020benders}}\gv{Rahmaniani et al. \cite{rahmaniani2020benders}} applied an inner approximation heuristic.
In the next section, we propose an alternative approach for generating Lagrangian cuts more quickly.

\section{Restricted Separation of Lagrangian Cuts}\label{sect:alg}
Observe that given any $(\pi,\pi_0)\in\bR^n\times\bR_+$, evaluating $\bQs(\pi,\pi_0)$ gives a valid inequality \eqref{cut}. By picking $(\pi,\pi_0)$ properly, we show in the following proposition that a small number of Lagrangian cuts is sufficient to give the perfect information bound \pv{\citep{avriel1970value}}\gv{\cite{avriel1970value}} defined as\begin{displaymath}
z_{\text{PI}}:=\sum_{s\in S}p_s\min_{x,y}\Big\{c^\top x+(q^s)^\top y:(x,y)\in K^s \Big\}.
\end{displaymath}
Although the perfect information bound can computed by solving a separate subproblem for each scenario, it is often better than the LP relaxation
bound because each subproblem is a MIP.
\begin{prop}\label{pibd}
	The following equality holds:\begin{equation}\label{eq:PI}
	z_{\text{PI}}=\min_{x,\theta_s}\Big\{c^\top x+\sum_{s\in S}p_s\theta_s:c^\top x+\theta_s\geq\bQs(c,1),s\in S \Big\}.
	\end{equation}
\end{prop}
\pv{\proof{Proof.}}
\gv{\begin{proof}}
	We first show the ``$\leq$" direction. It follows by observing that
if $(x,\{\theta_s\}_{s\in S})$ satisfies $c^\top x+\theta_s\geq\bQs(c,1)$ for all $s\in S$, then
	\begin{equation}\label{ineq:PI}
	\begin{aligned}
	c^\top x+\sum_{s\in S}p_s\theta_s\geq &~c^\top x+\sum_{s\in S}p_s(-c^\top x+\bQs(c,1))\\
	=&\sum_{s\in S}p_s\bQs(c,1)	=\sum_{s\in S}p_s\min_{x,y}\bigl\{c^\top x+(q^s)^\top y:(x,y)\in K^s
	\bigr\}=z_{\text{PI}}.
	\end{aligned}
	\end{equation}
	The inequality \eqref{eq:PI} follows from the fact that the first inequality in \eqref{ineq:PI} must hold at equality for all optimal solutions of \eqref{eq:PI}. This is because if $c^\top x+\theta_s>\bQs(c,1)$ for some $s\in S$, we can obtain a feasible solution with a lower objective value by decreasing the value of $\theta_s$.
	\pv{\Halmos}
\pv{\endproof}
\gv{\end{proof}}
From this example, we see that Lagrangian cuts with $(\pi,\pi_0)=(c,1)$ can already yield the perfect information bound. One interpretation of this example is that useful Lagrangian cuts can be obtained even when searching in a heavily constrained set of coefficients,
which motivates our study of restricted separation of Lagrangian cuts.

\subsection{General Framework}\label{gen_fw}
We propose solving the following restricted version of \eqref{separation} to search for a Lagrangian cut that separates a solution
$(\hat{x},\hat{\theta}_s)$:\begin{equation}\label{rstr_sepr}
\max_{\pi,\pi_0}\{\bQs(\pi,\pi_0) - \pi^\top \hat{x}-\pi_0\hat{\theta}_s:(\pi,\pi_0)\in\Pi_s \},
\end{equation}
where $\Pi_s$ can be any convex compact subset of $\bR^{n}\times\bR_+$, which is not necessarily defined by a neighborhood of the origin.

While setting $\Pi_s$ to be a neighborhood of the origin would give an exact separation of Lagrangian cuts, our proposal
is to use a more restricted definition of $\Pi_s$ with the goal of obtaining a separation problem that can be solved
faster. As one might expect, the strength of the cut generated from \eqref{rstr_sepr} will heavily depend on the choice
of $\Pi_s$. We propose to define $\Pi_s$ to be the span of a small set
of vectors. Specifically, we require\begin{displaymath}
\pi=\sum_{k=1}^K \beta_k\pi^k
\end{displaymath}
for some $\beta\in\bR^K$, where $\{\pi^k \}_{k=1}^K$ are chosen vectors. If $\{\pi^k \}_{k=1}^K$ span $\bR^n$, then
this approach reduces to exact separation. However, for smaller $K$, we expect a trade-off between computation time and
strength of cuts. A natural possibility for the vectors $\{\pi^k\}_{k=1}^K$, which we explore in Section \ref{Gen Pis}, is to use the coefficients of Benders cuts derived from the LP relaxations.

\subsection{Solution of \eqref{rstr_sepr}}\label{sect:soln}
We present a basic cutting-plane algorithm for solving \eqref{rstr_sepr} in Algorithm \ref{alg:rstr}.  The classical cutting-plane method used in the algorithm can be replaced by other bundle methods, e.g., the level method \pv{\citep{lemarechal1995new}}\gv{\cite{lemarechal1995new}}.
Within the solution process of \eqref{rstr_sepr}, the function $\bQs$ is approximated by an upper bounding cutting-plane model
$\hQs:\bR^{n}\times\bR_+ \rightarrow \bR$ defined by
\begin{equation}\label{Qcutpl}
\hQs(\pi,\pi_0)=\min_{z,\theta_s^z}\{\pi^\top z+\pi_0 \theta_s^z:(z,\theta_s^z)\in \hat{E}^s \}
\end{equation}
for $(\pi,\pi_0) \in\bR^n\times\bR_+$,
where $\hat{E}^s$ is a finite subset of $E^s$. 
In our implementation, for each scenario $s$, all feasible first-stage solutions obtained from evaluating
$\bQs(\pi,\pi_0)$ and the corresponding feasible second-stage costs are added into $\hat{E}^s$. Specifically, each time
we solve \eqref{subIP} for some $(\pi,\pi_0)$, we collect all optimal and sub-optimal solutions $(x,y)$ from the solver
and store the first-stage solution $z=x$ and the corresponding feasible second-stage cost $\theta_s^z=(q^s)^\top y$. To 
guarantee finite convergence, we assume that at least one of the optimal solutions obtained when solving \eqref{subIP} is an extreme point of
$\conv(K^s)$.
Note that even though we may solve problem \eqref{rstr_sepr} associated with different $\hat{x},\hat{\theta}_s$ and
$\Pi_s$, the previous $\hQs$ remains a valid overestimate of $\bQs$. Therefore, whenever we need to solve \eqref{subIP}
with some new $\hat{x},\hat{\theta}_s$ and $\Pi_s$, we use the currently stored $\hat{E}^s$ to initialize $\hQs$. In
our implementation, to avoid unboundedness of $\hat{Q}^*_{s}$, $\hat{E}^s$ is initialized by the perfect information
solution $(z^s,(q^s)^\top y^s)$, where $(z^{s},y^s)\in\arg\min_{x,y}\{c^\top x+(q^s)^\top y: (x,y) \in K^s \}$.

\begin{algorithm}[H]\label{alg:rstr}
	\SetAlgoLined
	\KwIn{($\hat{x},\hat{\theta}_s$), $\Pi_s$, $\hat{E}^s$}
	\KwOut{$(\pi^*,\pi^*_0)$, $\hat{E}^s$}
	Initialize UB$\leftarrow+\infty$, LB$\leftarrow-\infty$\\
	\While{UB $> 0$ and UB-LB$\geq\delta$UB}{
		Solve $UB\leftarrow\max_{\pi,\pi_0}\{\hQs(\pi,\pi_0)-\pi^\top \hat{x}-\pi_0\hat{\theta}_s:(\pi,\pi_0)\in\Pi_s\}$, and collect solution $(\pi,\pi_0)=(\hat{\pi},\hat{\pi}_0)$\\
		Solve \eqref{subIP} to evaluate $\bQs(\hat{\pi},\hat{\pi}_0)$. Update $\hat{E}^s$ and $\hQs$ with optimal and
		sub-optimal solutions obtained while solving \eqref{subIP}  \\
		\If{$LB<\bQs(\hat{\pi},\hat{\pi}_0)-\hat{\pi}^\top \hat{x}-\hat{\pi}_0\hat{\theta}_s$}
		{$LB\leftarrow\bQs(\hat{\pi},\hat{\pi}_0)-\hat{\pi}^\top \hat{x}-\hat{\pi}_0\hat{\theta}_s$\\
		$(\pi^*,\pi_0^*)\leftarrow(\hat{\pi},\hat{\pi}_0)$}
	}
	\caption{Solution of \eqref{rstr_sepr}.}
\end{algorithm}

While Algorithm \ref{alg:rstr} is for the most part a standard cutting-plane algorithm, an important detail  is the
stopping condition which uses a relative tolerance $\delta\in[0,1)$. Convergence of Algorithm \ref{alg:rstr} follows
standard arguments for a cutting-plane algorithm.  Specifically, because the set of extreme-point optimal solutions from
\eqref{subIP} is finite, there will eventually be an iteration where the optimal solution
obtained when solving \eqref{subIP} is already contained in the set $\hat{E}^s$. When that occurs, after updating $LB$
it will hold that $UB=LB$. At that point it either holds that UB $\leq 0$ or UB-LB $=0 < \delta$UB and the algorithm
terminates.
If the optimal value of \eqref{rstr_sepr} is positive, then UB will be strictly positive throughout the algorithm. In
that case Algorithm \ref{alg:rstr} returns a violating Lagrangian cut with cut violation
LB$\in((1-\delta)\text{UB},\text{UB}]$. The value of $\delta$ controls the trade-off between the strength of the cut
and the running time. On the other hand, if the optimal value of \eqref{rstr_sepr} is non-positive, then we can terminate
as soon as UB$\leq 0$, since in this case we are assured a violated cut will not be found.

We next demonstrate that even when using $\delta>0$ it is possible to attain the full strength of the Lagrangian cuts
for a given set $\Pi_s$ when applying them repeatedly within a cutting-plane algorithm for solving a relaxation of
problem \eqref{Benders}. To make this result concrete, we state such a cutting-plane algorithm in Algorithm \ref{alg:cutpl}. At each iteration $t$, we approximate the true objective function by the cutting-plane model:\begin{displaymath}
c^\top x+\sum_{s\in S}p_s\hat{Q}_s^t(x)=c^\top x+\sum_{s\in S}p_s\min_{x,\theta_s}\{\theta_s:\bar{\pi}^\top x+\bar{\pi}_0\theta_s\geq \bar{\tau},\ \forall (\bar{\pi}_0,\bar{\pi},\bar{\tau})\in\Psi_s^t \}
\end{displaymath}
where $\Psi_s^t$ represents the Benders cuts and Lagrangian cuts associated with scenario $s$ that have been added to
the master relaxation up to iteration $t$. At iteration $t$, a master problem based on the current cutting-plane model is solved and generates candidate
solution $(x^t,\{\theta_s^t \}_{s\in S})$. For each scenario $s$, the restricted Lagrangian cut separation problem
\eqref{rstr_sepr} is solved using Algorithm \ref{alg:rstr} with $\Pi_s$ restricted to be a (typically low-dimensional) set $\Pi_s^t$. The cutting-plane
model is then updated by adding the violated Lagrangian cuts (if any) for each scenario.

\begin{algorithm}[hbt!]
	\SetAlgoLined
	\KwIn{$\{\hat{Q}_s^0\}_{s\in S}$}
	\KwOut{$\{\hat{Q}_s^t\}_{s\in S}$}
	Initialize $t\leftarrow 0$, cutfound $\leftarrow$ True\\
	\While{cutfound = True}{ \label{alg:benstart} 
		cutfound $\gets$ False\;
		Solve $\min_{x,\theta_s}\{c^\top x+\sum_{s\in S}p_s\theta_s:\theta_s\geq\hat{Q}^t_s(x), Ax \geq b \}$ to obtain solution $(x^t,\{\theta_s^t \}_{s\in S})$\\
		\For{$s \in S$}{ 
			Solve Benders subproblem \eqref{BendersSub} and obtain dual solution $\mu^{t,s}$\; 
			\If{$(\mu^{t,s})^\top T^s x^t + \theta_s^t < (\mu^{t,s})^\top h^s$}{
				Update $\hat{Q}_s^t$ using \eqref{BenCut} defined by $\mu=\mu^{t,s}$ to obtain $\hat{Q}_s^{t+1}$\;
				cutfound $\gets$ True\; 
			}
		} \label{alg:benend} 
		\If{coutfound = False}{
		\For{$s\in S$}{
			Generate $\Pi_s^t\subseteq\bR^n$\\
			Solve \eqref{rstr_sepr} using Algorithm \ref{alg:rstr} with $(\hat{x},\hat{\theta}_s)=(x^{t},\theta_s^{t})$ and $\Pi_s=\Pi_s^t$, and collect solution $(\pi^{t,s},\pi^{t,s}_0):=(\pi^*,\pi^*_0)$\\
			\If{$\hat{Q}_s^*(\pi^{t,s},\pi_0^{t,s}) - (\pi^{t,s})^\top x^t - \pi^{t,s} \theta_s^t > 0$}{
			Update $\hat{Q}_s^{t}$ using cut \eqref{cut} with $(\pi,\pi_0)=(\pi^{t,s},\pi^{t,s}_0)$ to obtain
			$\hat{Q}_s^{t+1}$\;
			cutfound $\gets$ True\;	
			}
		}
		}
		$t\leftarrow t+1$
	}
	
	\caption{\THE}\label{alg:cutpl}
\end{algorithm}

We show in the next theorem that if we search on a static compact $\Pi_s$ and keep adding Lagrangian cuts by solving
\eqref{rstr_sepr} to a $\delta$ relative tolerance with $\delta\in[0,1)$, the algorithm will converge to a solution that
satisfies all the restricted Lagrangian cuts.

\begin{thm}\label{thm:convergence}
In Algorithm \ref{alg:cutpl}, let $\Pi_s^t\equiv\Pi_s$ be a compact set independent of $t$,
and assume that at each iteration of Algorithm \ref{alg:cutpl}, for each scenario $s\in S$, a Lagrangian cut (if a
violated one exists) is generated by solving \eqref{rstr_sepr} to $\delta$ relative tolerance with $\delta\in[0,1)$. Then for each scenario $s$, every accumulation point $(\bar{x},\bar{\theta}_s)$ of the sequence $\{(x^t,\theta_s^t)\}_{t=1}^\infty$ satisfies $(\bar{x},\bar{\theta}_s)\in P_s(\Pi_s)$.
\end{thm}
\pv{\proof{Proof.}}
\gv{\begin{proof}}
	Note that $\bQs$ is a piecewise linear concave function. Therefore, $\bQs$ is Lipschitz continuous.
	Assume $\{(x^{i_t},\theta_s^{i_t})\}_{t=1}^\infty$ is any fixed subsequence of $\{(x^t,\theta_s^t)\}_{t=1}^\infty$ that converges to a point $(\bar{x},\bar{\theta}_s)$.
	Assume for contradiction that\begin{displaymath}
	\max_{\pi,\pi_0}\{\bQs(\pi,\pi_0)-\pi^\top \bar{x}-\pi_0\bar{\theta}_s:(\pi,\pi_0)\in\Pi_s \}=:\epsilon>0.
	\end{displaymath}
	Since $\Pi_s$ is bounded, there exists $\tau$ such that\begin{displaymath}
	\max_{\pi,\pi_0}\{|\pi^\top (x^{i_t}-\bar{x})+\pi_0(\theta_s^{i_t}-\bar{\theta}_s)|:(\pi,\pi_0)\in\Pi_s\}\leq\epsilon/2\text{ for all }t\geq \tau.
	\end{displaymath}
	This implies that for all $t\geq \tau$,\begin{align*}
		&\max_{\pi,\pi_0}\{Q_s^*(\pi,\pi_0)-\pi^Tx^{i_t}-\pi_0\theta_s^{i_t}:(\pi,\pi_0)\in\Pi_s\}\\
		=&\max_{\pi,\pi_0}\{Q_s^*(\pi,\pi_0)-\pi^T\bar{x}-\pi_0\bar{\theta}_s-\pi^T(x^{i_t}-\bar{x})-\pi_0(\theta_s^{i_t}-\bar{\theta}_s):(\pi,\pi_0)\in\Pi_s\}\\
		\geq&\max_{\pi,\pi_0}\{Q_s^*(\pi,\pi_0)-\pi^T\bar{x}-\pi_0\bar{\theta}_s-|\pi^T(x^{i_t}-\bar{x})+\pi_0(\theta_s^{i_t}-\bar{\theta}_s)|:(\pi,\pi_0)\in\Pi_s\}\\
		\geq&\max_{\pi,\pi_0}\{Q_s^*(\pi,\pi_0)-\pi^T\bar{x}-\pi_0\bar{\theta}_s:(\pi,\pi_0)\in\Pi_s\}\\
		&-\max_{\pi,\pi_0}\{|\pi^T(x^{i_t}-\bar{x})+\pi_0(\theta_s^{i_t}-\bar{\theta}_s)|:(\pi,\pi_0)\in\Pi_s\}\\
		\geq &\epsilon-\epsilon/2=\epsilon/2,
	\end{align*}
	where the second ``$\geq$" follows from $\max_x\{f(x)\}\geq \max_x\{f(x)+g(x)\}-\max_x\{g(x)\}$.
	Let $(\pi^{(i_t)},\pi_0^{(i_t)})$ denote the multiplier associated with the Lagrangian cut added at iteration $i_t$. For any $t'>t\geq \tau$,\begin{displaymath}
\bQs(\pi^{(i_t)},\pi_0^{(i_t)}) -(\pi^{(i_t)})^\top x^{i_{t'}}-\pi_0^{(i_t)}\theta_s^{i_{t'}} \leq 0
	\end{displaymath}
	since the Lagrangian cut associated with $(\pi^{(i_t)},\pi_0^{(i_t)})$ was added before iteration $i_{t'}$. On the other hand,\begin{multline*}
	\bQs(\pi^{(i_{t'})},\pi_0^{(i_{t'})})-(\pi^{(i_{t'})})^\top x^{i_{t'}}-\pi_0^{(i_{t'})}\theta_s^{i_{t'}}\\
	\geq(1-\delta)\max_{\pi,\pi_0}\{\bQs(\pi,\pi_0)-\pi^\top x^{i_{t'}}-\pi_0\theta_s^{i_{t'}}:(\pi,\pi_0)\in\Pi_s\}\geq (1-\delta)\epsilon/2.
	\end{multline*}
	Since the sequence $\{(x^{i_t},\theta_s^{i_t})\}_{t=1}^\infty$ is bounded and function $\bQs$ is Lipschitz continuous, there exists $\rho>0$ such that $\|(\pi^{(i_t)},\pi_0^{(i_t)})-(\pi^{(i_{t'})},\pi_0^{(i_{t'})})\|\geq \rho$ for all $t'>t\geq \tau$. This contradicts with the fact that $\Pi_s$ is bounded.\pv{\Halmos}
\pv{\endproof}
\gv{\end{proof}}

In practice, we allow $\Pi_s^t$ to change with iteration $t$, but Theorem \ref{thm:convergence} provides motivation for solving \eqref{rstr_sepr} to a $\delta$ relative tolerance.
In Section \ref{sec:delta}, we present results of experiments demonstrating that for at least one problem class, using
$\delta$ much larger than 0 can lead to significantly faster improvement in the lower bound obtained when using Lagrangian cuts.
When allowing $\Pi_s^t$ to change with iteration $t$, Theorem \ref{thm:convergence} no longer applies. Instead, in this
setting we follow common practice and terminate the cut-generation process after  no more cuts are found, or
when it is detected that the progress in increasing the lower bound has significantly slowed (see Section \ref{sec:bnc} for an
example stopping condition). 

\subsection{Choice of $\Pi_s$}\label{Gen Pis}

The primary difference between the exact separation \eqref{separation} and the proposed restricted separation problem
\eqref{rstr_sepr} is the restriction constraint $(\pi,\pi_0)\in\Pi_s$. Ideally, we would like to choose $\Pi_s$ 
so that	\eqref{rstr_sepr} is easy to solve, in particular avoiding evaluating the function $\bQs$ too many times,
	while also yielding a cut \eqref{cut} that has a similar strength compared with the Lagrangian cut corresponding to the optimal solution of \eqref{separation}.
Our proposal that aims to satisfy these properties is to define $\Pi_s$ to be the span of past Benders cuts with some
normalization. Specifically, we propose
\begin{equation}\label{old_nml}
\Pi_s=\Big\{(\pi,\pi_0):\exists \beta\in\bR^K\text{ s.t. }\pi=\sum_{k=1}^K\beta_k\pi^k,\alpha\pi_0+\|\pi \|_1\leq 1,\pi_0\geq 0 \Big\}.
\end{equation}
Here $K$ and $\alpha$ are both predetermined parameters, and $\{\pi^k\}_{k=1}^K$ are the coefficients of the last $K$ Benders cuts corresponding to scenario $s$. The choice of $\{\pi^k\}_{k=1}^K$ is important. The motivation for
this choice is that the combination of cut coefficients from the LP relaxation could give a good approximation of the tangent of the supporting hyperplanes of $Q_s$ at an (near-)optimal solution. The normalization used here is similar to the one proposed for the selection of Benders cuts \pv{by}\gv{in} \cite{fischetti2010note}. 

Different normalization leads to different cuts. Our second proposal is to use \begin{equation}\label{new_nml}
\Pi_s=\Big\{(\pi,\pi_0):\exists \beta\in\bR^{K}\text{ s.t. }\pi=\sum_{k=1}^K\beta_k\pi^k,\alpha\pi_0+\|\beta\|_1\leq 1,\pi_0\geq 0 \Big\},
\end{equation}
which replaces $\|\pi \|_1$ in \eqref{old_nml} by $\|\beta \|_1$.
This normalization may tend to lead to solutions with sparse $\beta$. An advantage
of this normalization is that it leads to a MIP approximation for optimally choosing $\pi_k$'s from a set of candidate
vectors, which we discuss next.

Let $V^s=\{v^k\}_{k\in I_s}$ be a given (potentially large) collection of coefficient vectors for a scenario $s \in S$.   
We construct a MIP model for selecting a subset of these vectors of size at most $K$, to then use in \eqref{new_nml}.
Specifically, we introduce a binary variable $z_k$ for each $k \in I_s$, where $z_k = 1$ indicates we select $v^k$ as
one of the $\pi^k$ vectors to use in \eqref{new_nml}.  Using these decision variables,
the problem of selecting a subset of at most $K$ vectors that maximizes the normalized violation can be formulated as the following mixed-integer nonlinear program:
\begin{align}
\max_{\beta_k,z_k,\pi,\pi_0}\ &\bQs(\pi,\pi_0)-\pi^\top \hat{x}-\pi_0\hat{\theta}_s\label{pik selection}\\
\text{s.t.\ }\ &\pi=\sum_{k\in I_s}\beta_k v^k,\label{pik con1}\\
&\alpha\pi_0+\sum_{k\in I_s}|\beta_k|\leq 1, \label{lamnorm} \\
&-z_k\leq \beta_k\leq z_k, &k\in I_s, \label{lamzero} \\
&\sum_{k\in I_s}z_k\leq K, \label{card} \\
&z_k\in \{0,1\}, &k\in I_s,\\
&\pi_0\geq 0.\label{pik con-1}
\end{align}
The objective and constraints \eqref{pik con1} and \eqref{lamnorm} match the formulation \eqref{new_nml}, whereas
constraint \eqref{card} limits the number of selected vectors to at most $K$, and \eqref{lamzero} enforces that if a
vector is not selected then the weight on that vector must be zero.
This problem is computationally challenging due to the implicit function $\bQs$. However, if we replace $\bQs$ by its current cutting-plane approximation $\hQs$, then the problem can be solved as a MIP. Specifically, the MIP approximation of \eqref{pik selection}-\eqref{pik con-1} is given by:\begin{equation}\label{MIPapprox}
\begin{aligned}
\max_{\beta_k,z_k,\pi,\pi_0,\tau}\ &\tau - \pi^\top \hat{x}-\pi_0\hat{\theta}_s \\
\text{s.t.\ }\ &\tau\leq \pi^\top z+\pi_0\theta_s^z, &&(z,\theta_s^z)\in \hat{E}^s,\\
&\eqref{pik con1}-\eqref{pik con-1}.
\end{aligned}
\end{equation}
After solving this approximation, we choose the vectors $\{\pi_k\}_{k=1}^{K'}$ to be the vectors $v_k$ with $z_k=1$. 
Given this basis, we then solve \eqref{rstr_sepr} with $\Pi_s$ defined in \eqref{new_nml} to generate Lagrangian
cuts. The basis size $K'$ in this case can be strictly smaller than $K$. Observe that the optimal value of
\eqref{MIPapprox} gives an upper bound of \eqref{pik selection}-\eqref{pik con-1}. Therefore, if the optimal objective
value of \eqref{MIPapprox} is small, this would indicate we would not be able to generate a highly violated Lagrangian cut for this
scenario, and hence we can skip this scenario for generating Lagrangian cuts.


\section{Computational Study}\label{sec:CompStudy}

We report our results from a computational study investigating different variants of our proposed methods for generating
Lagrangian cuts and also comparing these to existing approaches. In Section \ref{sec:root} we investigate the bound
closed over time by just adding cuts to the LP relaxation. In Section \ref{sec:bnc} we investigate the impact of using
our cut-generation strategy to solve our test instances to optimality within a simple branch-and-cut algorithm.

We conduct our study on three test problems. The first two test problems are standard SIP problems that can be found in
literature: the stochastic server location problem (SSLP, binary first-stage and mixed-integer second-stage) introduced
\pv{by}\gv{in} \cite{ntaimo2005million}, and the stochastic network interdiction problem (SNIP, binary first-stage and
continuous second-stage) \pv{by}\gv{in} \cite{pan2008minimizing}. We also create a variant of the stochastic server
location problem (SSLPV) that has mixed-binary first-stage
variables and continuous second-stage varaibles. Our test set includes 24 SSLP instances with 20-50 first-stage
variables, 2000-2100
second-stage variables, and 50-200 scenarios. The SSLPV test instances are the same as for SSLP, except the number of
first-stage variable is doubled, with half of them being continuous.
We use 20 instances of the SNIP problem having 320 first-stage variables, 2586 second-stage
variables, and 456 scenarios. More details of the test problems and instances can be found in \pv{Section S.\ref{test_prob} of the online supplement}\gv{Appendix}.

We consider in total six implementations of Lagrangian cuts:\begin{enumerate}
	\item {\strben}: Strengthened Benders cuts from \cite{zou2019stochastic}, i.e., \eqref{ineq:Lag} with $\lambda$ equal to the negative of the Benders cut coefficients each time a Benders cut is generated.
	\item {\BDD}: \BDD$_3$ from \cite{rahmaniani2020benders}. This algorithm uses a multiphase implementation that first generates Benders cuts and strengthened Benders cuts, and then generates Lagrangian cuts based on a regularized inner approximation for separation.
	\item {\exact}: Exact separation of Lagrangian cuts of the form \eqref{rstr_sepr} with $\Pi_s=\{(\pi,\pi_0): \alpha\pi_0+\|\pi\|_1\leq 1,\pi_0\geq 0 \}$.
	\item {\rstrone}: Restricted separation of Lagrangian cuts \eqref{rstr_sepr} with $\Pi_s$ defined in \eqref{old_nml} and $\{\pi_k\}_{k=1}^{K}$ being the coefficients of the last $K$ Benders cuts corresponding to scenario $s$. 
	\item {\rstrtwo}: Restricted separation of Lagrangian cuts \eqref{rstr_sepr} with $\Pi_s$ defined in \eqref{new_nml} and $\{\pi_k\}_{k=1}^{K}$ being the coefficients of the last $K$ Benders cuts corresponding to scenario $s$. 
	\item {\rstrmip}: Restricted separation of Lagrangian cuts \eqref{rstr_sepr} with $\Pi_s$ defined in \eqref{new_nml} and $\{\pi_k\}_{k=1}^{K}$ determined by solving the MIP approximation \eqref{MIPapprox} with $V^s=\{v^k\}_{k\in I_s}$ being the set of all Benders cuts coefficients that have been generated for scenario $s$.
\end{enumerate}

\subsection{Implementation Details}\label{sec:impl}

All experiments are run on a Windows laptop with 16GB RAM and an Intel Core i7-7660U processor running at 2.5GHz. All
LPs, MIPs and convex quadratic programs (QPs) are solved using the optimization solver Gurobi 8.1.1 for SNIP and SSLP instances, and are solved using Gurobi 9.0.3 for SSLPV instances.

In Algorithm \ref{alg:rstr}, we also terminate if the upper bound is small enough ($<10^{-6}(|\hat{\theta}_s|+1)$) or two consecutive multipliers are too close to each other (infinity norm of the difference $<10^{-10}$).
A computational issue that can occur when solving \eqref{rstr_sepr} is that if we obtain a solution with $\pi_0=0$ (or very
small), the solution of $(x^*,y^*)$ of \eqref{subIP} does not provide information
about $Q_s(x^*)$ since the coefficients associated with the $y$ variables are equal zero. In that case, if we store the
``optimal" second-stage cost $(q^s)^\top y^*$ and use it together with $x^*$ to define a cut in the cutting-plane model
\eqref{Qcutpl}, this second-stage cost value $(q^s)^\top y^*$ can be very far away from $Q_s(x^*)$ and lead to a very loose
cut. So in our implementation, when $\pi_0$ is small ($<10^{-4}$), after finding $(x^*,y^*)$ from solving \eqref{subIP}, we reevaluate $Q_s(x^*)$ by solving the scenario subproblem \eqref{def:Qs} with $x=x^*$ fixed and use this value as $\theta_s^{x^*}$.

For adding Benders cuts in Algorithm \ref{alg:cutpl}, the initialization of $\hat{Q}_s^0$ is obtained by one iteration
of Benders cuts. Specifically, to avoid unboundedness, we solve $\min_x\{0:Ax\geq b\}$ to obtain a candidate solution
and add the corresponding Benders cuts for all scenarios. 
As described in Algorithm \ref{alg:cutpl}, lines \ref{alg:benstart}-\ref{alg:benend}, 
Benders cuts from the LP relaxation are added using the classical cutting-plane method
\pv{\citep{kelley1960cutting}}\gv{\cite{kelley1960cutting}}. This implementation was used for the SNIP test problem.
For the SSLP test problem we found this phase took too long to converge, so the phase of adding Benders cuts as
described in lines \ref{alg:benstart}-\ref{alg:benend}
was replaced by the level method \pv{\citep{lemarechal1995new}}\gv{\cite{lemarechal1995new}}. 
Benders cuts are only added if the violation is at least $10^{-4}(|\hat{\theta}_s|+1)$. We store the
coefficients of all Benders cuts that have been added to the master problem. For each $s\in S$, we define $\Pi_s$ in the
definitions of {\rstrone} and {\rstrtwo} based on the coefficients of the $K$ Benders cuts that have been most recently
found by the algorithm for scenario $s$.

We use Algorithm \ref{alg:rstr} for solving the separation problem, with one exception. When doing exact separation of
Lagrangian cuts (i.e., without constraining $\Pi_s$ to be low-dimensional set) we implement a level method for solving
\eqref{rstr_sepr}, because in that case we found Algorithm \ref{alg:rstr} converged too slowly. Based on preliminary
experiments, the normalization coefficient $\alpha$ in \eqref{old_nml} or \eqref{new_nml} is set to be 1 for tests on SSLP and SSLPV instances, and 0.1 for tests on SNIP instances. 

In all SSLP, SSLPV, and SNIP instances the set $\{x:Ax\leq b\}$ is integral and we have relatively complete recourse. Therefore, no feasibility cuts are necessary. Only Lagrangian cuts with $\pi_0\geq 10^{-6}$ are added in our tests.

For the {\BDD} tests, we follow the implementation of \cite{rahmaniani2020benders} and set
the regularization coefficient as $\delta_0=0.01$ for SSLP and SSLPV instances and as $\delta_0=100$ for SNIP instances. These
values are chosen based on preliminary experiments as producing the best improvement in bound over time.

We have also tested our instances using the general purpose SIP solver DSP
\pv{\citep{kim2018algorithmic,kim2019asynchronous,dandurandscalable}}\gv{\cite{kim2018algorithmic,kim2019asynchronous,dandurandscalable}}.
All DSP experiments are run on an Ubuntu virtual machine built on the same Windows machine using the binary version of
DSP 1.1.3 with CPLEX 12.10 as the MIP solver.

\subsection{LP Relaxation Results}\label{sec:root}

We first compare results obtained from different implementations of Lagrangian cuts on our test problems and
study the impact of different parameters in the algorithm. The stopping condition of Algorithm \ref{alg:cutpl} for all
root node experiments is that either the algorithm hits the one-hour time limit or there is no cut that can be found by the algorithm (with a $10^{-6}$ relative tolerance).

To compare the overall performance of the methods visually across many test instances, we present results in the form of
a $\gamma$-gap-closed profile. Given a set of problem instances $P$ and a set of cut generation methods $M$, let
$\bar{g}_p$ denote the largest gap closed by any of the methods (compared with the LP bound) in $M$ for instance $p\in
P$. The $\gamma$-gap-closed profile is defined with respect to certain threshold $\gamma\in(0,1]$. Given the value of
$\gamma$, we define $t^\gamma_{p,m}$ as the earliest time of closing the gap by at least $\gamma\bar{g}_p$ for problem $p\in P$ with method $m\in M$. If method $m$ did not close a gap at least $\gamma\bar{g}_p$ before termination, we define $t^\gamma_{p,m}=\infty$. The $\gamma$-gap-closed profile is a figure representing the cumulative growth (distribution function) of the $\gamma$-gap-closed ratio $\rho_m^\gamma(\tau)$'s over time $\tau$ where\begin{displaymath}
\rho_m^\gamma(\tau)=\frac{|\{p\in P: t^\gamma_{p,m}\leq \tau \}|}{|P|}.
\end{displaymath}
In this paper, $P$ is either the set of all SSLP test instances, all SNIP test instances or all SSLPV instances, and $M$ is the set of all
methods with all parameters we have tested. Therefore, for each instance $p$, $\bar{g}_p$ is defined to be the largest
gap we have been able to close with Lagrangian cuts. We set parameter $\gamma$ equal to either $0.75$ or $0.95$ for
summarizing the information regarding time needed for obtaining what we refer to as a ``reasonably good bound''
($\gamma=0.75$) or a ``strong bound'' ($\gamma=0.95$).

To illustrate the behavior of these algorithms we also present the lower bound obtained over time using different cut
generation schemes on some individual instances. 

\subsubsection{Relative Gap Tolerance $\delta$\pv{.}}\label{sec:delta}

\begin{figure}[hbt!]
	\centering
	\begin{subfigure}[b]{0.46\linewidth}
		\pv{\includegraphics[width=\linewidth]{Exact_delta_SSLP.png}}
		\gv{\includegraphics[width=\linewidth]{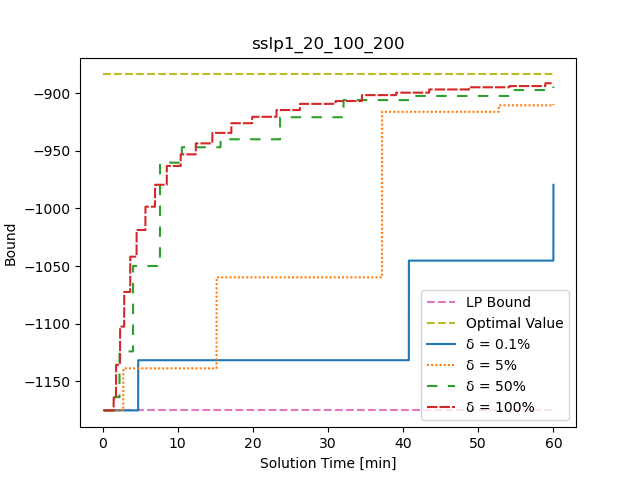}}
	\end{subfigure}
	\begin{subfigure}[b]{0.46\linewidth}
		\pv{\includegraphics[width=\linewidth]{Exact_delta_SSLPV.png}}
		\gv{\includegraphics[width=\linewidth]{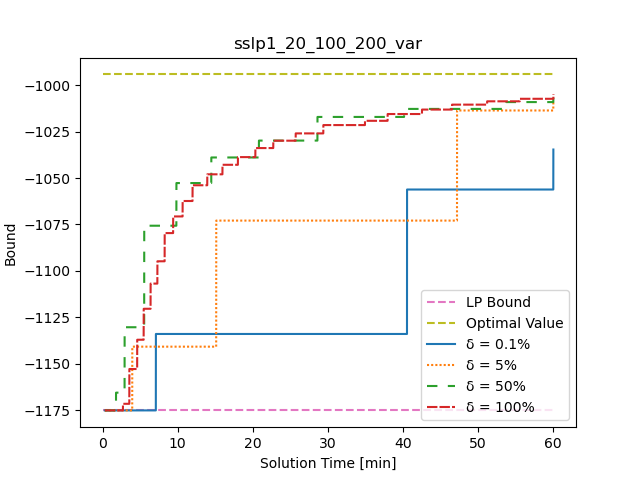}}
	\end{subfigure}
	\caption{Convergence profile of {\exact} on an SSLP instance (left) {and an SSLPV instance (right)} with varying
	$\delta$ values.}
	\label{fig:delta_Exact}
\end{figure}
\begin{figure}[hbt!]
	\centering
	\begin{subfigure}[b]{0.46\linewidth}
		\pv{\includegraphics[width=\linewidth]{SSLP_Exact_delta_75.png}}
		\gv{\includegraphics[width=\linewidth]{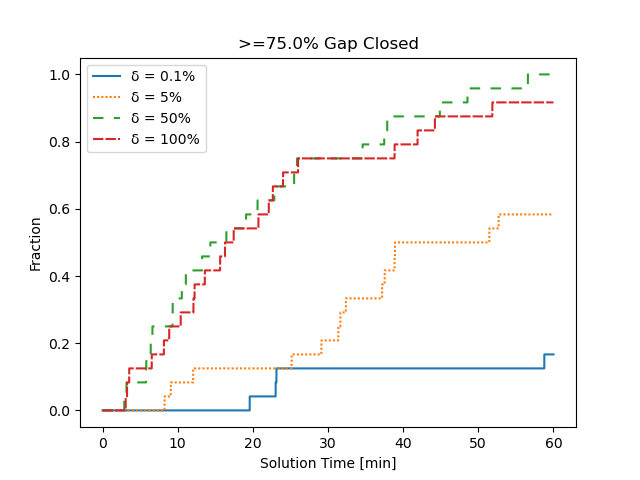}}
	\end{subfigure}
	\begin{subfigure}[b]{0.46\linewidth}
		\pv{\includegraphics[width=\linewidth]{SSLP_Exact_delta_95.png}}
		\gv{\includegraphics[width=\linewidth]{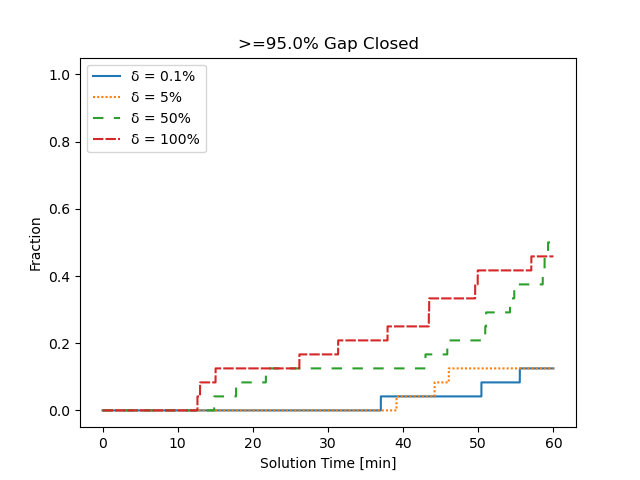}}
	\end{subfigure}
	\begin{subfigure}[b]{0.46\linewidth}
		\pv{\includegraphics[width=\linewidth]{SSLP_Rstr1_delta_75.png}}
		\gv{\includegraphics[width=\linewidth]{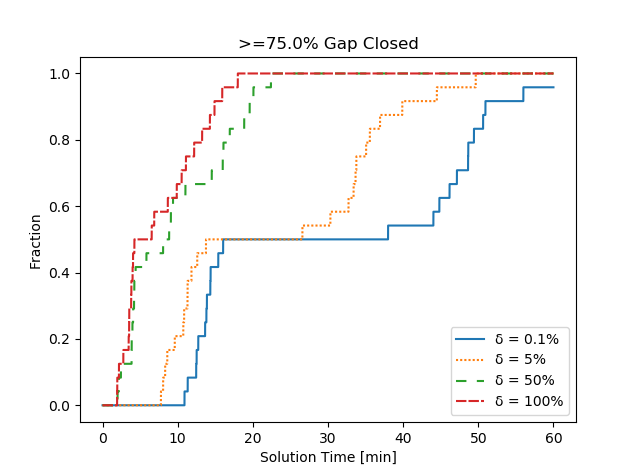}}
	\end{subfigure}
	\begin{subfigure}[b]{0.46\linewidth}
		\pv{\includegraphics[width=\linewidth]{SSLP_Rstr1_delta_95.png}}
		\gv{\includegraphics[width=\linewidth]{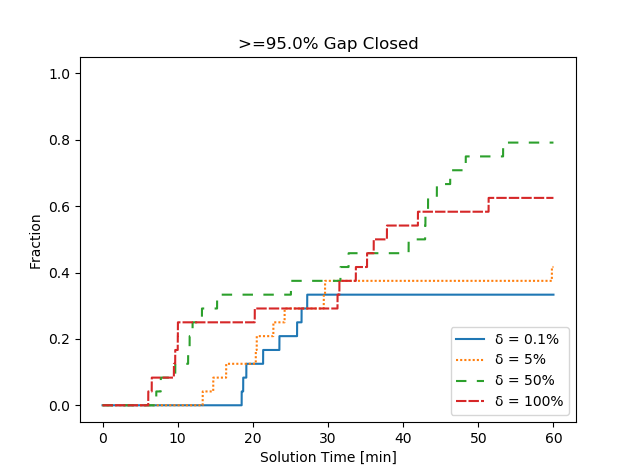}}
	\end{subfigure}
	\caption{$\gamma$-gap-closed profile for SSLP instances obtained by {\exact} (top) and {\rstrone} ($K=10$, bottom)
	with $\gamma=0.75$ (left) or $\gamma=0.95$ (right) and varying $\delta$ values.}
	\label{fig:SSLP_delta}
\end{figure}
\begin{figure}[hbt!]
	\centering
	\begin{subfigure}[b]{0.46\linewidth}
		\pv{\includegraphics[width=\linewidth]{SNIP_Rstr1_delta_75.png}}
		\gv{\includegraphics[width=\linewidth]{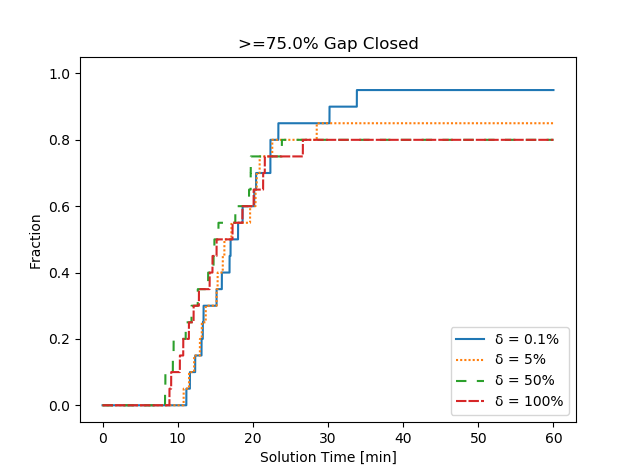}}
	\end{subfigure}
	\begin{subfigure}[b]{0.46\linewidth}
		\pv{\includegraphics[width=\linewidth]{SNIP_Rstr1_delta_95.png}}
		\gv{\includegraphics[width=\linewidth]{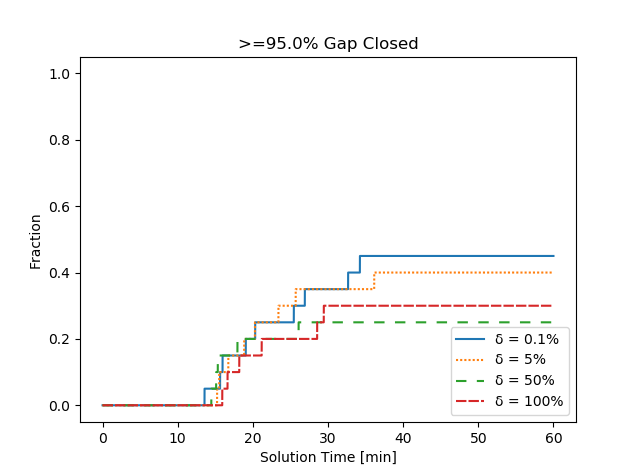}}
	\end{subfigure}
	\caption{$\gamma$-gap-closed profile for SNIP instances obtained by {\rstrone} with $\gamma=0.75$ (left) or
	$\gamma=0.95$ (right), $K=10$ and varying $\delta$ values.}
	\label{fig:SNIP_delta}
\end{figure}
We first examine the impact of varying the relative gap tolerance $\delta$. Experiments are run with parameter $\delta$ being $0.1\%, 5\%, 50\%$ and $100\%$, respectively. A small $\delta$ value (e.g.,
$0.1\%$) means solving \eqref{rstr_sepr} to near optimality, which may generate strong cuts but can be time-consuming to
generate each cut. Using large $\delta$ values might sacrifice some strength of the cut in an iteration for an earlier termination.

For different $\delta$ values, in Figure \ref{fig:delta_Exact}, we compare the bounds obtained by {\exact} over time. As
the trend is quite similar for instances within the same problem class, we only report the convergence profiles for one
SSLP instance (sslp1\_20\_100\_200) and one SSLPV instance (sslp1\_20\_100\_200\_var). For the SNIP instances, which have many
more first-stage variables, {\exact} often cannot finish one iteration of separation within an hour, so we put the convergence profile of {\exact} for one SNIP instance in \pv{Section S.\ref{cvg_profiles} of the online supplement}\gv{Appendix}. For the SSLP and SSLPV
instances, a large $\delta$ value tends to close the gap faster than a small $\delta$ value. Plots providing the same
comparison for {\rstrone} with $K=10$ are provided in \pv{Section S.\ref{cvg_profiles} of the online supplement}\gv{Appendix}. For {\rstrone}, a
large $\delta$ value helps improve the bound efficiently for SSLP and SSLPV instances while the choice of $\delta$ has little
impact on the performance of the SNIP instances.

To summarize the impact of $\delta$ over the full set of test instances, we plot the $0.75$-gap-closed profile and $0.95$-gap-closed profile with different $\delta$ values for SSLP and SNIP instances in Figures
\ref{fig:SSLP_delta} and \ref{fig:SNIP_delta}. We omit the gap-closed profiles for the SSLPV instances since they are
similar to those for the SSLP instances. The gap-closed profiles for SNIP instances with {\exact} are omitted as
no useful bound is obtained by {\exact} for any of the SNIP instances. For the SSLP and SSLPV instances, it is preferable to use a
large $\delta$ value as both a reasonable bound ($0.75$ gap closed) and a strong bound ($0.95$ gap closed) are obtained earlier when $\delta$ is large. For the
SNIP instances, there is less clear preference although small $\delta$ values perform slightly better. This is because
the time spent on each iteration of separation does not vary too much by changing the $\delta$ value for the SNIP instances. 

\noindent\textbf{Conclusion.} Depending on the problem, the choice of $\delta$ can have either a mild or big impact. But we find a relatively large $\delta$, e.g., $\delta=50\%$, has a stable performance on our instance classes.

\subsubsection{Basis Size $K$\pv{.}}\label{sec:K}
\begin{figure}[hbt!]
	\centering
	\begin{subfigure}[b]{0.46\linewidth}
		\pv{\includegraphics[width=\linewidth]{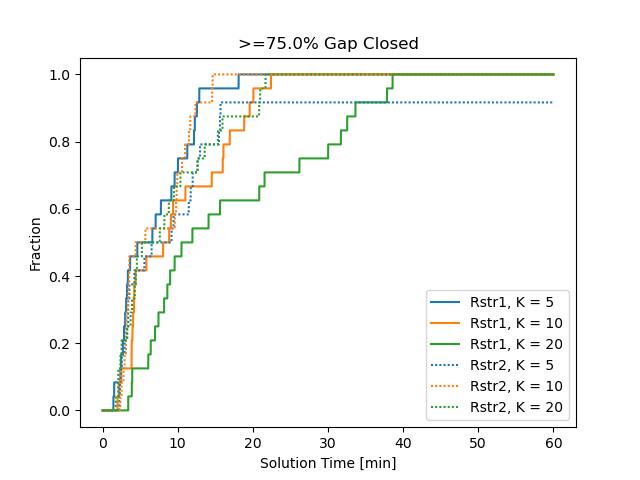}}
		\gv{\includegraphics[width=\linewidth]{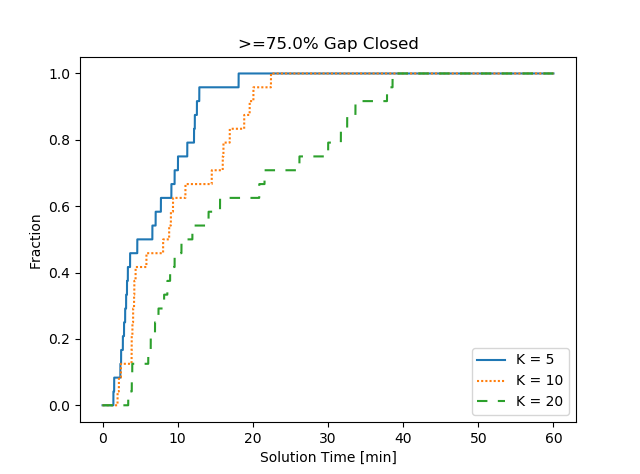}}
	\end{subfigure}
	\begin{subfigure}[b]{0.46\linewidth}
		\pv{\includegraphics[width=\linewidth]{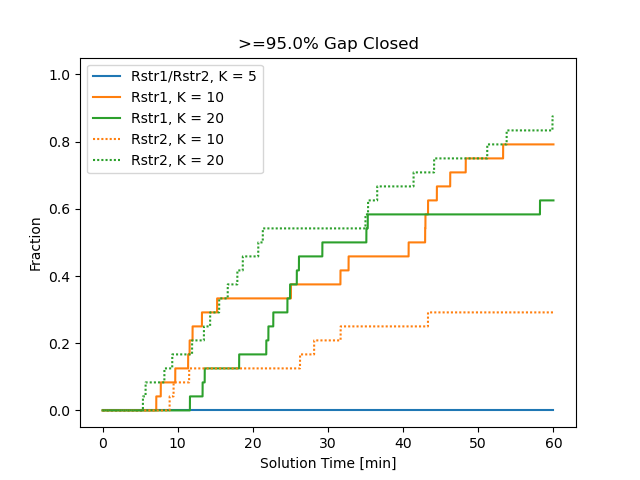}}
		\gv{\includegraphics[width=\linewidth]{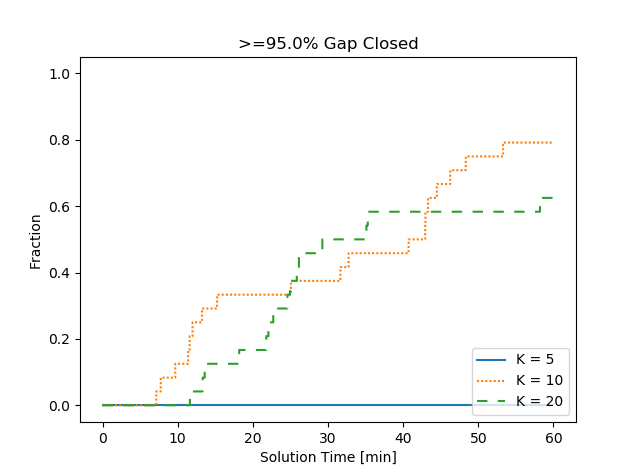}}
	\end{subfigure}
	\begin{subfigure}[b]{0.46\linewidth}
		\pv{\includegraphics[width=\linewidth]{SSLP_RstrMIP_K_75.png}}
		\gv{\includegraphics[width=\linewidth]{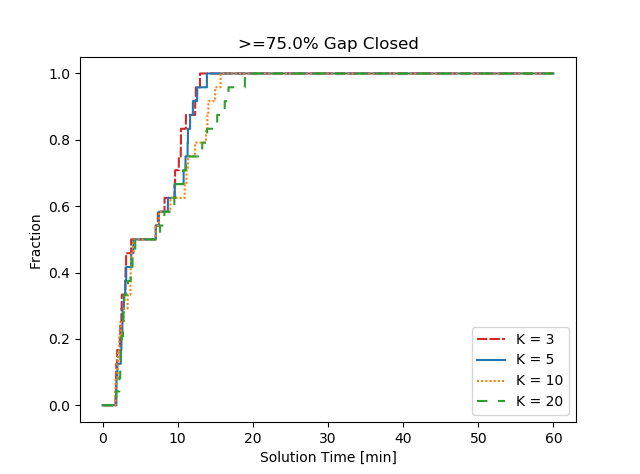}}
	\end{subfigure}
	\begin{subfigure}[b]{0.46\linewidth}
		\pv{\includegraphics[width=\linewidth]{SSLP_RstrMIP_K_95.png}}
		\gv{\includegraphics[width=\linewidth]{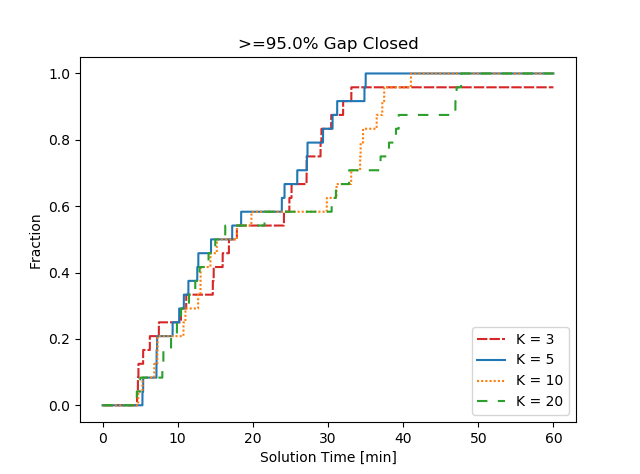}}
	\end{subfigure}
	\caption{$\gamma$-gap-closed profile for SSLP instances obtained by {\rstrone} and {\rstrtwo} (top), and {\rstrmip}
	(bottom) with $\gamma=0.75$ (left) or $\gamma=0.95$ (right), $\delta=50\%$ and varying $K$ values.}
	\label{fig:SSLP_K}
\end{figure}
\begin{figure}[hbt!]
	\centering
	\begin{subfigure}[b]{0.46\linewidth}
		\pv{\includegraphics[width=\linewidth]{SNIP_Rstr12_K_75.png}}
		\gv{\includegraphics[width=\linewidth]{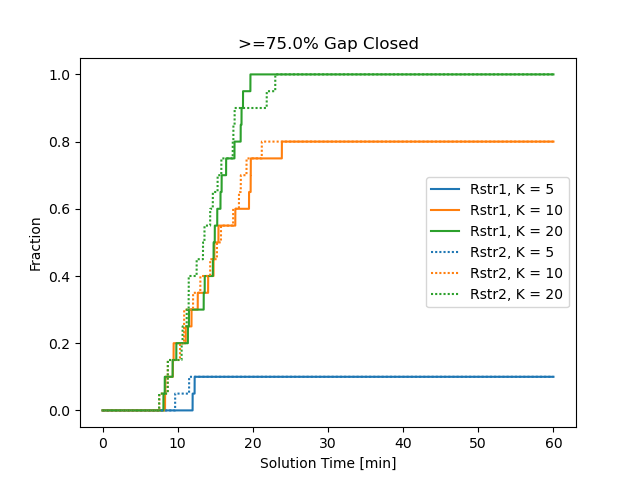}}
	\end{subfigure}
	\begin{subfigure}[b]{0.46\linewidth}
		\pv{\includegraphics[width=\linewidth]{SNIP_Rstr12_K_95.png}}
		\gv{\includegraphics[width=\linewidth]{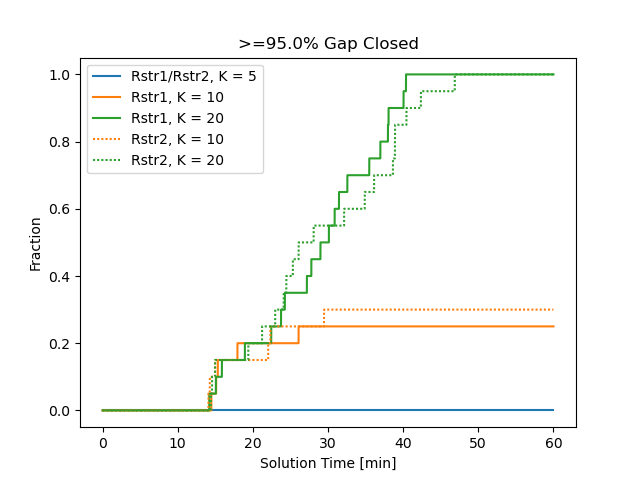}}
	\end{subfigure}
	\begin{subfigure}[b]{0.46\linewidth}
		\pv{\includegraphics[width=\linewidth]{SNIP_RstrMIP_K_75.png}}
		\gv{\includegraphics[width=\linewidth]{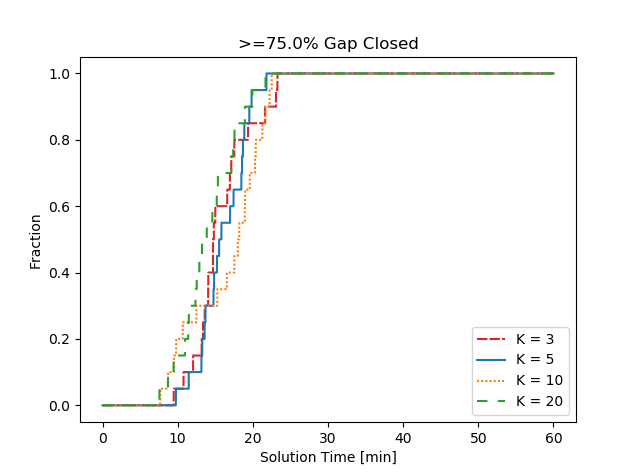}}
	\end{subfigure}
	\begin{subfigure}[b]{0.46\linewidth}
		\pv{\includegraphics[width=\linewidth]{SNIP_RstrMIP_K_95.png}}
		\gv{\includegraphics[width=\linewidth]{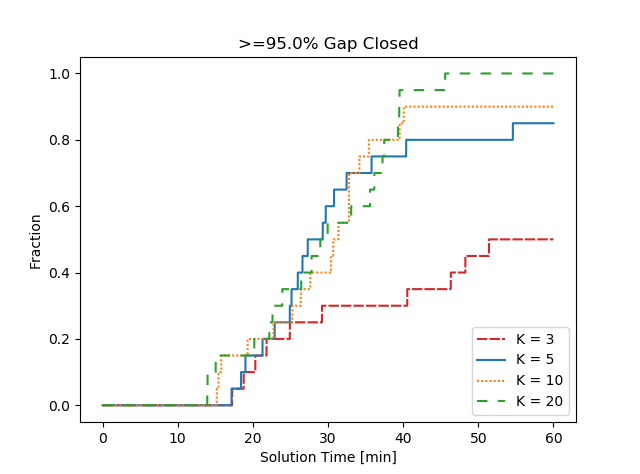}}
	\end{subfigure}
	\caption{$\gamma$-gap-closed profile for SNIP instances obtained by {\rstrone} and {\rstrtwo} (top),  and {\rstrmip}
	(bottom) with $\gamma=0.75$ (left) or $\gamma=0.95$ (right), $\delta=50\%$ and varying $K$ values.}
	\label{fig:SNIP_K}
\end{figure}
We next investigate the impact of the choice of the basis size $K$ in {\rstrone}, {\rstrtwo} and {\rstrmip}. In Figures
\ref{fig:SSLP_K} and \ref{fig:SNIP_K}  we plot the $0.75$-gap-closed profiles and the $0.95$-gap-closed profiles with
different $K$ values for SSLP and SNIP instances, respectively. The results for SSLPV instances are omitted as they are similar to the results for SSLP instances. Plots showing the evolution of the lower bound
obtained using {\rstrone}, {\rstrtwo} and {\rstrmip} with different $K$ values on an example SNIP instance, an
example SSLP instance and an example SSLPV instance are given in \pv{Section S.\ref{cvg_profiles} of the online supplement}\gv{Appendix}.

As seen in Figures \ref{fig:SSLP_K} and \ref{fig:SNIP_K}, the trends are quite different for different problem classes when applying {\rstrone}. For SSLP and SSLPV instances, smaller $K$ performs better for obtaining a reasonably good
bound faster since the restricted separation is much easier to solve when $K$ is small. However, for SNIP instances
$K=20$ seems to dominate $K=5$ and $K=10$. We believe this is because the time to generate a cut is not significantly
impacted by $K$ for the SNIP instances, but $K=20$ gives much stronger cuts for closing the gap. 
Comparing {\rstrone} and {\rstrtwo}, {\rstrtwo} appears to be less sensitive to $K$ in terms of obtaining a reasonably good bound for SSLP instances. On the other hand, the performance of {\rstrtwo} is almost identical to the performance of {\rstrone} for SNIP instances. For SSLP and SSLPV instances with {\rstrone} or {\rstrtwo}, we find that larger $K$ often leads to longer solution time for obtaining a reasonable bound while small $K$ ($K=5$) is incapable of getting a strong bound.
We observe that {\rstrmip} is fairly robust to the choice of the basis size $K$. One particular reason is that the constraint $\sum_{k\in I_s}z_k\leq K$ in \eqref{MIPapprox} is often loose in early stages of the algorithm for relatively large $K$'s as the constraint $\alpha\pi_0+\|\beta\|_1\leq 1$ tends to select a sparse $\beta$ and thus a sparse $z$. For the same instance, different basis sizes result in very close bounds at the end of the run.

\noindent\textbf{Conclusion.} The basis size $K$ yields a trade-off between the strength of cuts and computation time. We find the optimal choice of $K$ for {\rstrone} depends on the instance. Slightly larger $K$ ($K=20$) often works better for {\rstrtwo} while the performance of {\rstrmip} is much more robust to the choice of $K$.

\subsubsection{Comparison between Methods\pv{.}}\label{Compare}
\begin{figure}[hbt!]
	\centering
	\begin{subfigure}[b]{0.46\linewidth}
		\pv{\includegraphics[width=\linewidth]{SSLP75.png}}
		\gv{\includegraphics[width=\linewidth]{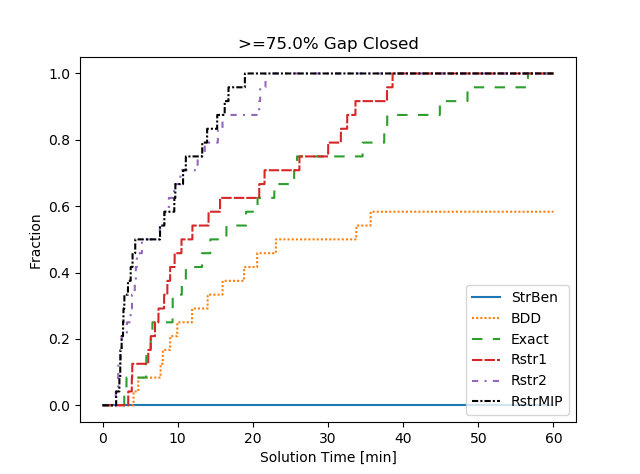}}
	\end{subfigure}
	\begin{subfigure}[b]{0.46\linewidth}
		\pv{\includegraphics[width=\linewidth]{SSLP95.png}}
		\gv{\includegraphics[width=\linewidth]{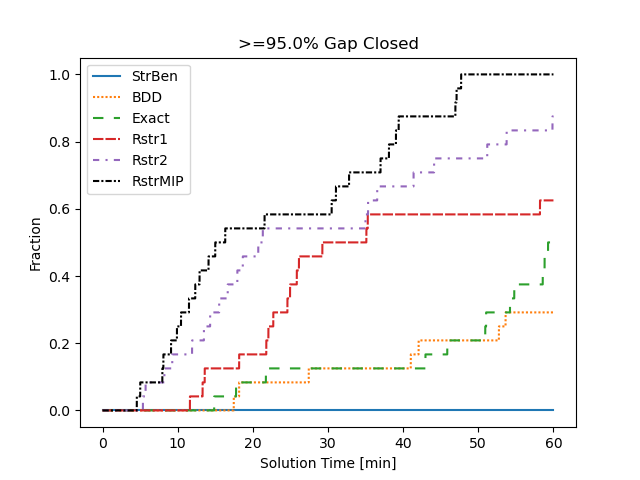}}
	\end{subfigure}
\begin{subfigure}[b]{0.46\linewidth}
		\pv{\includegraphics[width=\linewidth]{SSLPV75.png}}
		\gv{\includegraphics[width=\linewidth]{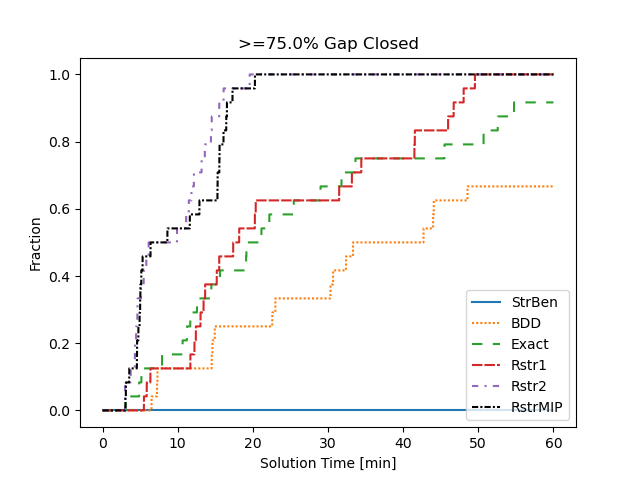}}
	\end{subfigure}
			\begin{subfigure}[b]{0.46\linewidth}
		\pv{\includegraphics[width=\linewidth]{SSLPV95.png}}
		\gv{\includegraphics[width=\linewidth]{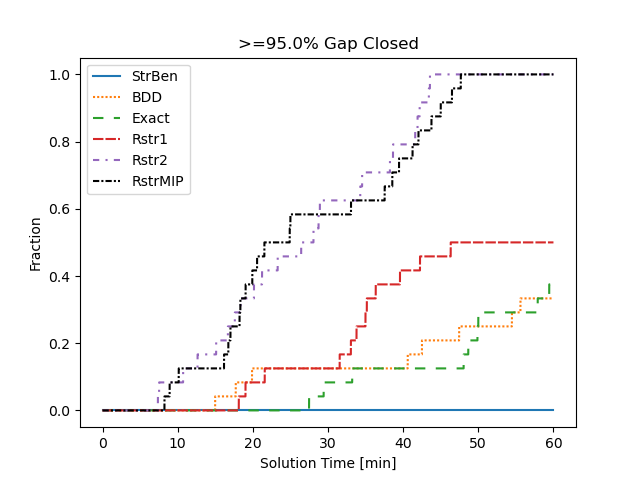}}
	\end{subfigure}
\begin{subfigure}[b]{0.46\linewidth}
		\pv{\includegraphics[width=\linewidth]{SNIP75.png}}
		\gv{\includegraphics[width=\linewidth]{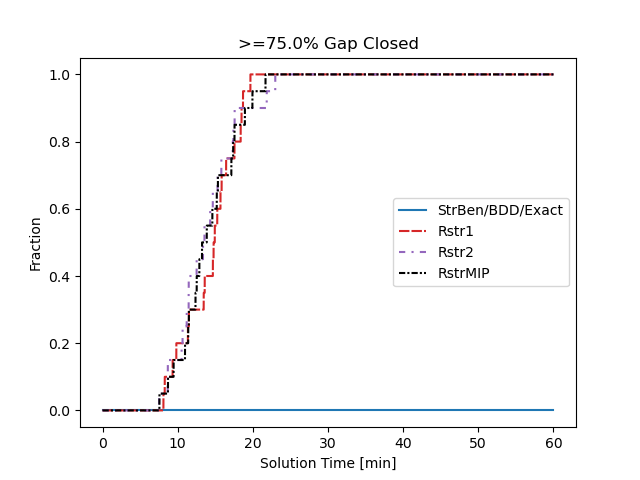}}
	\end{subfigure}
	\begin{subfigure}[b]{0.46\linewidth}
		\pv{\includegraphics[width=\linewidth]{SNIP95.png}}
		\gv{\includegraphics[width=\linewidth]{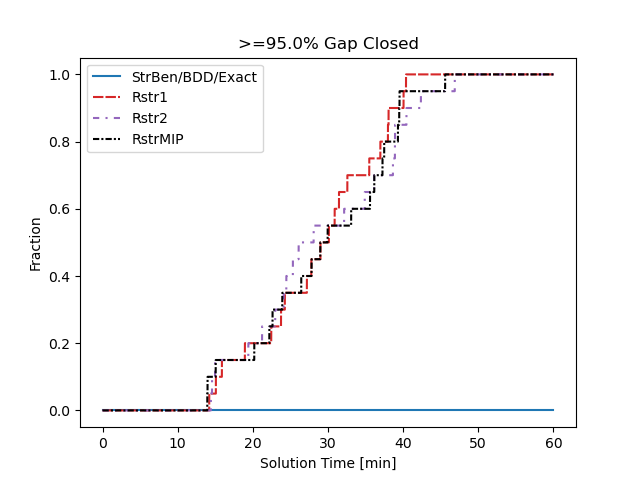}}
	\end{subfigure}
	\caption{$\gamma$-gap-closed profile for SSLP instances (top), {SSLPV instances (middle)}, and SNIP instances (bottom)
	with $\gamma=0.75$ (left) or $\gamma=0.95$ (right).}
	\label{fig:GapProf}
\end{figure}

We present results obtained by strengthened Benders cut ({\strben}), the Benders dual decomposition ({\BDD}) and
different variants of our restricted Lagrangian cut separation ({\exact}, {\rstrone}, {\rstrtwo} and {\rstrmip}). For
all variants of restricted Lagrangian cut separation, we choose relative gap tolerance $\delta=50\%$. For basis sizes,
we set $K=20$ for {\rstrone}, {\rstrtwo} and {\rstrmip}. We plot $0.75$-gap-closed profiles and $0.95$-gap-closed
profiles for these methods in Figure \ref{fig:GapProf}.

Apparently, {\rstrmip} has the best performance on SSLP instances in terms of the gap closed over time. 
On SSLPV instances, the performances of {\rstrtwo} and {\rstrmip} are
similar while much better than other methods. 
On SNIP
instances, the performance of {\rstrone}, {\rstrtwo} and {\rstrmip} are almost identical while {\rstrmip} is much more
robust in the choice of $K$ as observed in Section \ref{sec:K}. The curves for {\BDD} are flat on the bottom plots of
Figure \ref{fig:GapProf} as {\BDD} fails to close 75\% of the gap for all SNIP instances. Convergence profiles for
{\BDD} on one SSLP instance, one SSLPV instance and one SNIP instance can be found in \pv{Section S.\ref{cvg_profiles} of the online supplement}\gv{Appendix}.

\noindent\textbf{Conclusion.} {\rstrmip} has the best overall performance and in particular is more robust to the choice of $K$.

\subsubsection{Integrality gaps}
In Table \ref{Table:IntGapClosed}, we report statistics (minimum, maximum and average) of the LP integrality gap and
the integrality gap obtained after 60 minutes of Lagrangian cut generation by {\rstrmip} with $K=20$ for each problem class. The integrality gaps
are calculated as ($z_{\text{IP}}-$LB)/$|z_{\text{IP}}|$, where LB is the lower bound obtained from the Benders model
before or after adding the generated Lagrangian cuts. The results show that {\rstrmip} closes the vast majority of the
gap to the optimal value. Since the resulting bound closely approximates the optimal
value, and the exact Lagrangian dual must lie between this bound and the optimal value, these results imply that, on
these test instances, our method approximates the optimal value of the Lagrangian dual, and that the Lagrangian dual approximates the optimal
value. We emphasize, however, that this experiment is just to illustrate the potential for this
method to close a majority of the gap when given enough time -- we recommend using an early stopping condition to stop the cut generation process earlier when it is detected
that the bound improvement has slowed significantly.

\pv{\begingroup\renewcommand{\arraystretch}{0.67}}
\begin{table}[hbt!]
	\begin{center}
		\caption{Integrality gap closed by Lagrangian cuts generated by {\rstrmip}.}\label{Table:IntGapClosed}
		\begin{tabular}{ @{\extracolsep{\fill}}
		lp{0.08\textwidth}p{0.08\textwidth}p{0.09\textwidth} p{0.02\textwidth}  p{0.08\textwidth}p{0.08\textwidth}p{0.09\textwidth}}
			\toprule
			 & & & & & \multicolumn{3}{c}{Gap after {\rstrmip}}\\ 
			& \multicolumn{3}{c}{LP gap (\%)} & & \multicolumn{3}{c}{ cuts (\%)}\\
			Class & min & max & average & & min & max & average\\
			\cmidrule(lr){2-4}\cmidrule(lr){6-8}\addlinespace
			SSLP & 32.9 & 58.0 & 44.4 & & 0.0 & 3.6 & 0.7\\
			SSLPV & 18.2 & 27.8 & 23.8 & & 0.0 & 2.1 & 0.8\\
			SNIP & 22.1 & 34.2 & 28.0 & & 0.8 & 5.9 & 3.0\\
			\bottomrule
		\end{tabular}
	\end{center}
\end{table}
\pv{\endgroup}

\subsection{Tests for Solving to Optimality}\label{sec:bnc}
Finally, we study the use of our approach for generating Lagrangian cuts within a branch-and-cut method  to solve SIPs to
optimality, where Benders cuts and integer L-shaped cuts are added using Gurobi lazy constraint callback when a new
solution satisfying the integrality constraints is found. We refer to this method as LBC. We consider only {\rstrmip}
with $\delta=50\%$ and $K=10$ as the implementation of Lagrangian cuts. To reduce the time spent on Lagrangian cut
generation at the root node, an early termination rule is applied. We stop generating Lagrangian cuts if the gap closed
in the last 5 iterations is less than 1\% of the total gap closed so far. For comparison, we also experiment with
Benders cuts integrated with branch-and-cut (BBC -- \pv{algorithm described in Section S.\ref{BnC} of the online supplement}\gv{Algorithm \ref{alg:BnC} in Appendix}) and the DSP solver
\pv{\citep{kim2018algorithmic}}\gv{\cite{kim2018algorithmic}} on the same instances. A two-hour time limit is set for
these experiments. Gurobi heuristics are turned off for the branch-and-cut master problem to avoid the bug in Gurobi
8.1.1 that the MIP solution generated by Gurobi heuristics sometimes cannot trigger the callback. (This bug is supposed
to be fixed in Gurobi 9.1.)
\begin{table}[hbt!]
	\begin{center}
		\caption{Comparison of algorithms for solving SSLP instances to optimality.}\label{Table:SSLP}
		\begin{tabular}{ @{\extracolsep{\fill}} lrrrrrr}
			\toprule
			& \multicolumn{3}{c}{On $|S|=50$ instances} & \multicolumn{3}{c}{On $|S|=200$ instances}\\
			& LBC & BBC & DSP & LBC & BBC & DSP\\
			\cmidrule(lr){2-4}\cmidrule(lr){5-7}\addlinespace
			\# solved & 12/12 & 0/12 & 12/12 & 12/12 & 0/12 & 7/12\\
			Avg soln time & 1496 & $\geq 7200$ & 1259 & 4348 & $\geq 7200$ & $\geq 5365$\\
			Avg gap (\%) & 0.0 & 18.1 & --- & 0.0 & 27.0 & ---\\
			Avg B\&C time & 25 & $\geq 7200$ & --- & 163 & $\geq 7200$ & ---\\
			Avg \# nodes & 775 & 13808 & --- & 1411 & 3954 & ---\\
			\bottomrule
		\end{tabular}
	\end{center}
\end{table}
\pv{\endgroup}

 We present the results obtained by LBC, BBC and DSP on SSLP instances in Table \ref{Table:SSLP}. We report the number of instances at sample size $|S|=50$ or $|S|=200$ each algorithm solved within the time
limit (`\# solved'), the total average solution time of each algorithm (`Avg soln time') including time to generate cuts if relevant, the average ending optimality gap (`Avg gap (\%)'), average time spent solving the problem with branch-and-cut {\it not including the initial cut generation time} (`Avg B\&C time'), and the average number of branch-and-bound nodes explored (`Avg \#
nodes').
The ending optimality gap of a method on an instance is calculated as
$(\text{UB}-\text{LB})/\max\{|\text{UB}|,|\text{LB}|\}$ where UB and LB are the upper and lower bounds obtained from the method when the time limit was reached.

All SSLP instances are solved to optimality by LBC within the two-hour time limit with at most a few thousand nodes
explored. On the other hand, none of the SSLP instances are solved by BBC before hitting the time limit, leaving a relatively large optimality gap. SSLP instances are hard to solve by BBC since neither Benders cuts or integer L-shaped cuts are strong
enough to cut off infeasible points efficiently in the case when second-stage variables are discrete. Comparing with
DSP, LBC solves more $|S|=200$ SSLP instances than DSP withinin the time limit. The performance of LBC is relatively stable compared with DSP in the
sense that LBC scales well as $|S|$ increases. The efficiency of DSP heavily depends on the strength of the
nonanticipative dual bound. For all instances where DSP outperforms LBC on solution time, no branching is needed
for closing the dual gap. We observe that LBC spends most of the solution time on generating
Lagrangian cuts for SSLP instances. The branch-and-cut time required after generating Lagrangian cuts is very small.

\pv{\begingroup\renewcommand{\arraystretch}{0.67}}
\begin{table}[hbt!]
	\begin{center}
		\caption{Comparison of algorithms for solving SSLPV instances to optimality.}\label{Table:SSLPV}
		\begin{tabular}{ @{\extracolsep{\fill}} lrrrrrr}
			\toprule
			& \multicolumn{3}{c}{On $|S|=50$ instances} & \multicolumn{3}{c}{On $|S|=200$ instances}\\
			& LBC & BBC & DSP & LBC & BBC & DSP\\
			\cmidrule(lr){2-4}\cmidrule(lr){5-7}\addlinespace
			\# solved & 11/12 & 0/12 & 3/12 & 7/12 & 0/12 & 2/12\\
			Avg soln time & $\geq 2537$ & $\geq 7200$ & --- & $\geq 5894$ & $\geq 7200$ & ---\\
			Avg gap (\%) & 0.1 & 18.4 & --- & 0.3 & 20.8 & ---\\
			Avg B\&C time & $\geq 605$ & $\geq 7197$ & --- & $\geq 1156$ & $\geq 7191$ & ---\\
			Avg \# nodes & $\geq 7101$ & $\geq 22970$ & --- & $\geq 6300$ & $\geq 8945$ & ---\\
			\bottomrule
		\end{tabular}
	\end{center}
\end{table}
\pv{\endgroup}

In Table \ref{Table:SSLPV}, we report results obtained by these methods for the SSLPV instances. We do not report the
average solution time of DSP because for many of the SSLPV instances DSP fails at a very early stage by outputting an
infeasible solution. Even though the performances of our approaches  at the root node  is similar for SSLPV as SSLP, 
we find the SSLPV instances are far more difficult to solve to optimality than the SSLP instances.
As with the SSLP instances, BBC is unable to solve any of the SSLPV instances because the Benders cuts are not strong enough.
Within the two-hour time limit, only a small fraction of the SSLPV instances can be solved by DSP correctly. On the
other hand, LBC solves the majority of the SSLPV instances, and the ending optimality gaps are small for those unsolved
instances, showing its capability and stability in dealing with different problem classes.

\pv{\begingroup\renewcommand{\arraystretch}{0.67}}

\begin{table}[hbt!]
	\begin{center}
		\caption{Comparison of algorithms for solving SNIP instances to optimality.}\label{Table:SNIP}
		\begin{tabular}{ @{\extracolsep{\fill}} lrrr}
			\toprule
			& LBC & BBC & DSP\\
			\hline
			\# solved & 20/20 & 20/20 & 0/20\\
			Avg soln time & 2333 & 163 & $\geq 7200$\\
			Avg gap (\%) & 0.0 & 0.0 & ---\\
			Avg B\&C time & 71 & 116 & ---\\
			Avg \# nodes & 1231 & 3348 & ---\\
			\bottomrule
		\end{tabular}
	\end{center}
\end{table}


\begin{table}[hbt!]
	\begin{center}
		\caption{Comparison of LBC and BBC on SNIP instances  with Gurobi presolve and cuts turned off.}\label{BnCbasic:SNIP}
		\begin{tabular}{ @{\extracolsep{\fill}} lrr}
			\toprule
			& LBC & BBC\\
			\hline
			\# solved & 20/20 & 9/20 \\
			Avg soln time & 2385 & $\geq 4981$\\
			Avg gap (\%) & 0.0 & 2.9\\
			Avg B\&C time & 71 & $\geq 4534$\\
			Avg \# nodes & 4236 & $\geq 2923914$\\
			\bottomrule
		\end{tabular}
	\end{center}
\end{table}
\pv{\endgroup}

In Table \ref{Table:SNIP} we report results using these methods to solve the SNIP instances. We observe that no SNIP instance can be solved within two hours by DSP. This is because solving the nonanticipative dual
problem of SNIP instances is not completed in two hours. BBC solves the SNIP instances faster than LBC as the time spent
on cut generation at the root node is much shorter and the Benders cuts, when combined with Gurobi's presolve and cut
generation procedures are sufficient for keeping the number of nodes small.

To better understand the impact of the Lagrangian cuts vs.~Benders cuts on their own, we also report results obtained by
LBC and BBC with Gurobi presolve and cuts disabled in Table \ref{BnCbasic:SNIP}. We see that Gurobi presolve and cuts are quite essential for solving
SNIP instances when using only Benders cuts. In this case, BBC is unable to solve many of the instances within the time limit. This is because Gurobi presolve and
cuts significantly improve the LP bounds for BBC. LBC still solves all SNIP instances, demonstrating the strength of the Lagrangian cuts generated. 

\section{Conclusion}\label{sec:conclusion}
We derive a new method for generating Lagrangian cuts to add to a Benders formulation of a two-stage SIP. 
Extensive numerical experiments on three SIP problem classes indicate that the proposed method can improve the bound 
more quickly than alternative methods for using Lagrangian cuts.
When incorporating the proposed method within a simple branch-and-cut algorithm for solving the problems to optimality
we find that our method is able to solve one problem class within a time limit that is
is not solved by the open-source solver DSP, and performs modestly better on the other problem classes.

A direction of future work is to explore the integration of our techniques for generating Lagrangian cuts with the use
of other classes of cuts for two-stage SIPs, and also for multi-stage SIPs via the SDDIP algorithm 
\pv{proposed by  \cite{zou2019stochastic}}\gv{proposed by \cite{zou2019stochastic}}
